# Discovering New Runge–Kutta Methods Using Unstructured Numerical Search

By

David K. Zhang

Thesis

Submitted to the Faculty of Vanderbilt University

in partial fulfillment of the requirements for the degree of

BACHELOR OF SCIENCE

with Honors in Mathematics

April 16, 2019

Nashville, Tennessee

| Approved: | Date: |
|---|---|
| Douglas Hardin, Ph.D. | |
| Robert Womersley, Ph.D. | |
| Wöden Kusner, Ph.D. | |

**Abstract**

Runge–Kutta methods are a popular class of numerical methods for solving ordinary differential equations. Every Runge–Kutta method is characterized by two basic parameters: its *order*, which measures the accuracy of the solution it produces, and its *number of stages*, which measures the amount of computational work it requires. The primary goal in constructing Runge–Kutta methods is to maximize order using a minimum number of stages. However, high-order Runge–Kutta methods are difficult to construct because their parameters must satisfy an exponentially large system of polynomial equations. This paper presents the first known $10^{\text{th}}$-order Runge–Kutta method with only 16 stages, breaking a 40-year standing record for the number of stages required to achieve $10^{\text{th}}$-order accuracy. It also discusses the tools and techniques that enabled the discovery of this method using a straightforward numerical search.

*This thesis is dedicated to my mother and father, whose love, like all good mathematical objects, exists and is unique.*

# TABLE OF CONTENTS

Chapter 1

# Runge–Kutta Methods

## 1.1 A Historical Introduction

Runge–Kutta methods are a popular class of numerical methods for solving initial value problems for ordinary differential equations. They are widely employed by computational scientists in physics [2, 5, 8, 46], astronomy [22, 47], chemistry [37, 39, 45], biology [35, 49, 55], pharmacology [1], economics [7, 59], and even in non-scientific fields, including video game development [27, 56]. Thus, the construction of high-quality Runge–Kutta methods is a field of considerable practical importance, which has been studied for over 120 years [54].

As with all numerical approximation algorithms, the primary goal in developing Runge–Kutta methods is to produce highly accurate results using a small amount of computational effort. In the language of the Runge–Kutta literature, we aim to construct methods of high *order* using few *stages*. (Precise definitions for these terms will be given in Section 1.2.) This is a difficult balance to strike because the parameters of a high-order Runge–Kutta method need to satisfy an exponentially large number of constraints. In particular, a method with $s \in \mathbb{N}$ stages has $s(s+1)/2$ free parameters, but to achieve order $p \in \mathbb{N}$, these parameters must satisfy a system of approximately $3^p$ polynomial equations, known as *order conditions*. The huge gap between quadratic and exponential growth causes the system of order conditions to quickly become overdetermined for any practically reasonable number of stages.

In addition, these order conditions exhibit a complex recursive structure that makes them notoriously difficult to analyze by hand. As a famous example, in 1965 [13], John Butcher proved a lower bound on the number of stages required for a Runge–Kutta method to achieve order $p = 7$. It then took Butcher another twenty years to complete the proof of the corresponding lower bound for order $p = 8$, which was finally published in 1985 [15]. (These bounds will be presented as Theorem 1.27 in Section 1.4 of this paper; however, their proofs are far too complicated to be reproduced here.)

Historically, the standard approach to overcoming these algebraic difficulties has been "the skillful use of simplifying assumptions" to reduce the number of equations that must be simultaneously considered [34]. That is, instead of directly solving the full system of $\approx 3^p$ order conditions, one constructs a *reduced system* containing fewer equations that collectively form a sufficient condition for the full system. This reduced system is designed to be simple enough to solve by hand, but restrictive enough to imply all necessary order conditions. Of course, this approach may prevent certain methods from being discovered if the reduced system forms a sufficient but not necessary condition. The construction of these reduced systems typically involves more art than science, often with essential properties determined by arbitrary assumptions and unfounded heuristics. (See Section II.5 of [34] for several horrifying examples.)

In 1975, A. R. Curtis [20] became a leader in this enterprise by constructing the first[1] Runge–Kutta method of order 10, achieved using 18 stages. Three years later, in 1978, Ernst Hairer stole the lead "by using the complete arsenal of simplifying ideas" to construct a method of order 10 using only 17 stages. This set the record for the highest-order Runge–Kutta method constructed in the 20$^{\text{th}}$ century. Although methods of higher order have been found in the 21$^{\text{st}}$ century [24, 25], no methods of order 10 have ever been constructed using fewer than 17 stages.

In this paper, we improve upon Hairer's result by presenting the first Runge–Kutta method of order 10 that uses only 16 stages. We describe the process by which this method was discovered *numerically*, without complicated manual analysis or simplifying assumptions of any kind. We discuss the variety of algorithmic tools and numerical techniques that have enabled this discovery on modern computer hardware, most of which (to the best of our knowledge) have never previously been applied to the discovery of Runge–Kutta methods.

We note that this is not the first time a novel Runge–Kutta method has been discovered by numerical means. In 2009, Sergey Khashin obtained the first 13-stage Runge–Kutta method of order 9 by applying Tikhonov-regularized Newton iteration to a reduced system of order conditions [41]. However, we emphasize that Khashin's work relied upon a system of "filtrated simplifying assumptions" that he introduced to reduce the dimensionality of the problem. In contrast, our work demonstrates the feasibility of directly solving the full system of order conditions, eliminating the possibility that simplifying assumptions may prevent some methods from being discovered.

The rest of the paper is organized as follows. The remainder of Chapter 1 provides an expository account of the theory of Runge–Kutta methods and Butcher's derivation of the system of order conditions using rooted trees. Chapter 2 describes our strategy for numerically solving the resulting system of polynomial equations, along with our tactics for accelerating the search using the facilities of modern computer processors. Finally, Chapter 3 compares the method that was produced by this procedure to previous Runge–Kutta methods and provides concluding remarks.

[1]Technically, this should be called the first *nontrivial* Runge–Kutta method of order 10, since the order of any Runge–Kutta method can be arbitrarily increased by Richardson extrapolation. However, this procedure produces methods with an excessive number of stages. For example, applying 10$^{\text{th}}$-order Richardson extrapolation to Euler's method produces a Runge–Kutta method with 46 stages.

## 1.2 Preliminary Definitions

**Definition 1.1.** For the purposes of this paper, an ***initial value problem*** in dimension $n \in \mathbb{N}$ is an ordered triple $(\mathbf{f}, t_0, \mathbf{y}_0)$ consisting of a smooth function $\mathbf{f} : \mathbb{R} \times \mathbb{R}^n \to \mathbb{R}^n$ and two initial points, $t_0 \in \mathbb{R}$ and $\mathbf{y}_0 \in \mathbb{R}^n$. A ***solution*** of an initial value problem is a smooth function $\mathbf{y} : \mathbb{R} \to \mathbb{R}^n$ satisfying

$$\mathbf{y}'(t) = \mathbf{f}(t, \mathbf{y}(t)) \qquad \text{and} \qquad \mathbf{y}(t_0) = \mathbf{y}_0. \tag{1.1}$$

An ***autonomous*** initial value problem is one in which $\mathbf{f}$ does not depend on $t$. In this case, we write $\mathbf{f} : \mathbb{R}^n \to \mathbb{R}^n$ and require that the solution $\mathbf{y} : \mathbb{R}^n \to \mathbb{R}^n$ satisfy

$$\mathbf{y}'(t) = \mathbf{f}(\mathbf{y}(t)) \qquad \text{and} \qquad \mathbf{y}(t_0) = \mathbf{y}_0. \tag{1.2}$$

We will often speak of the independent variable $t$ as representing time, and we will interpret $\mathbf{f}$ as prescribing the evolution of the dependent variable $\mathbf{y}$ over time. This is purely a matter of convention adopted for linguistic convenience, and our analysis of initial value problems will apply equally well to problems with other types of independent variables.

Smoothness of $\mathbf{f}$ and $\mathbf{y}$ on all of $\mathbb{R} \times \mathbb{R}^n$ is an unnecessarily restrictive hypothesis, but we will make this standing assumption throughout this paper because it does not affect the development of Runge–Kutta methods. Of course, to apply these methods effectively in practice, we only need local differentiability of $\mathbf{f}$ in a neighborhood of $(t_0, \mathbf{y}_0)$. In particular, for a Runge–Kutta method to achieve convergence of order $k$, we need $\mathbf{f}$ to be locally of class $C^{k+1}$ [10].

**Remark 1.2.** Any non-autonomous initial value problem $(\mathbf{f}, t_0, \mathbf{y}_0)$ in dimension $n$ can be rewritten as an autonomous initial value problem in dimension $n + 1$ by defining $\mathbf{z}_0 \in \mathbb{R}^{n+1}$ and $\mathbf{g} : \mathbb{R}^{n+1} \to \mathbb{R}^{n+1}$ as follows:

$$\mathbf{z}_0 := \begin{bmatrix} t_0 \\ \mathbf{y}_0 \end{bmatrix} \qquad \mathbf{g}(\mathbf{z}) := \begin{bmatrix} 1 \\ \mathbf{f}(z_1, \mathbf{z}_{2:n+1}) \end{bmatrix}$$

Here, $\mathbf{z}_{2:n+1}$ denotes the vector in $\mathbb{R}^n$ obtained from $\mathbf{z} \in \mathbb{R}^{n+1}$ by dropping its first component. The idea is to fold the time parameter $t$ into the state vector $\mathbf{z}$ and prescribe its time-evolution as $\mathrm{d}t/\mathrm{d}t = 1$. A solution of $(\mathbf{g}, t_0, \mathbf{z}_0)$ then corresponds to a solution of $(\mathbf{f}, t_0, \mathbf{y}_0)$ by dropping its first component. We call $(\mathbf{g}, t_0, \mathbf{z}_0)$ the ***autonomous form*** of $(\mathbf{f}, t_0, \mathbf{y}_0)$.

Runge–Kutta methods are a class of algorithms that take as input an initial value problem $(\mathbf{f}, t_0, \mathbf{y}_0)$ and a step size parameter $h \in \mathbb{R}$. They evaluate $\mathbf{f}$ a fixed number of times on linear combinations of $t_0$, $h$, $\mathbf{y}_0$, and previous evaluations of $\mathbf{f}$, and output an approximation $\mathbf{y}_1$ of the exact solution at a later time, $\mathbf{y}_1 \approx \mathbf{y}(t_0 + h)$. By repeatedly applying a Runge–Kutta method to the new triple $(\mathbf{f}, t_0 + h, \mathbf{y}_1)$, we can advance the approximate solution returned by a Runge–Kutta method to arbitrarily many time points.

The structure of Runge–Kutta methods is best illustrated by example. The simplest possible Runge–Kutta method is ***Euler's method***, which we will present in the following slightly unusual notation:

$$\mathbf{k}_1 := h\,\mathbf{f}(t_0, \mathbf{y}_0) \tag{1.3}$$

$$\mathbf{y}_1 := \mathbf{y}_0 + \mathbf{k}_1 \tag{1.4}$$

(In what follows, we will always denote by $\mathbf{y}_1$ the approximate solution produced by a Runge–Kutta method.) There are several refinements of Euler's method that give more accurate solutions for a slight increase in computational work. These include the ***midpoint method*** and ***Heun's method***, presented below on the left and right, respectively.

$$\mathbf{k}_1 := h\,\mathbf{f}(t_0, \mathbf{y}_0) \tag{1.5}$$

$$\mathbf{k}_2 := h\,\mathbf{f}(t_0 + \tfrac{1}{2}h, \mathbf{y}_0 + \tfrac{1}{2}\mathbf{k}_1) \tag{1.6}$$

$$\mathbf{y}_1 := \mathbf{y}_0 + \mathbf{k}_2 \tag{1.7}$$

$$\mathbf{k}_1 := h\,\mathbf{f}(t_0, \mathbf{y}_0) \tag{1.8}$$

$$\mathbf{k}_2 := h\,\mathbf{f}(t_0 + h, \mathbf{y}_0 + \mathbf{k}_1) \tag{1.9}$$

$$\mathbf{y}_1 := \mathbf{y}_0 + \tfrac{1}{2}\mathbf{k}_1 + \tfrac{1}{2}\mathbf{k}_2 \tag{1.10}$$

The ***classical Runge–Kutta method***, often abbreviated as RK4, is the most famous member of this class. In fact, it is often referred to as "*the* Runge–Kutta method" without further elaboration [10].

$$\mathbf{k}_1 := h\mathbf{f}(t_0, \mathbf{y}_0) \tag{1.11}$$

$$\mathbf{k}_2 := h\mathbf{f}(t_0 + \frac{1}{2}h, \mathbf{y}_0 + \frac{1}{2}\mathbf{k}_1) \tag{1.12}$$

$$\mathbf{k}_3 := h\mathbf{f}(t_0 + \frac{1}{2}h, \mathbf{y}_0 + \frac{1}{2}\mathbf{k}_2) \tag{1.13}$$

$$\mathbf{k}_4 := h\mathbf{f}(t_0 + h, \mathbf{y}_0 + \mathbf{k}_3) \tag{1.14}$$

$$\mathbf{y}_1 := \mathbf{y}_0 + \frac{1}{6}\mathbf{k}_1 + \frac{1}{3}\mathbf{k}_2 + \frac{1}{3}\mathbf{k}_3 + \frac{1}{6}\mathbf{k}_4 \tag{1.15}$$

As illustrated by these examples, Runge–Kutta methods all exhibit the same logical structure: a sequence of evaluations of $\mathbf{f}$ to obtain a fixed number of intermediate stages $\{\mathbf{k}_i\}$, followed by a linear combination of these stages to obtain $\mathbf{y}_1$. Thus, a Runge–Kutta method is specified by a number of stages together with the coefficients employed at each stage. The standard notation for these coefficients is codified in the following definition.

**Definition 1.3.** Let $s \in \mathbb{N}$. An *s-**stage (explicit)***[2] ***Runge–Kutta method*** consists of a strictly lower-triangular matrix $A \in \mathbb{R}^{s\times s}$ and two vectors $\mathbf{b}, \mathbf{c} \in \mathbb{R}^s$, and is applied to an initial value problem $(\mathbf{f}, t_0, \mathbf{y}_0)$ with step size $h \in \mathbb{R}$ as follows:

$$\mathbf{k}_i := h\mathbf{f}\left(t_0 + c_i h, \mathbf{y}_0 + \sum_{j=1}^{i-1} a_{ij}\mathbf{k}_j\right) \qquad \text{for} \qquad i = 1, 2, \dots, s \tag{1.16}$$

$$\mathbf{y}_1 := \mathbf{y}_0 + \sum_{i=1}^{s} b_i \mathbf{k}_i \tag{1.17}$$

The parameters $(A, \mathbf{b}, \mathbf{c})$ of a Runge–Kutta method are traditionally displayed in a ***Butcher table*** of the following form:

$$\begin{array}{c|c} \mathbf{c} & A \\ \hline & \mathbf{b}^T \end{array} \tag{1.18}$$

To illustrate this notational convention, Butcher tables for the Euler (1.3)–(1.4), midpoint (1.5)–(1.7), Heun (1.8)–(1.10), and classical Runge–Kutta (1.11)–(1.15)

[2] The Runge–Kutta methods discussed in this paper are called "explicit methods" in the literature to distinguish them from "implicit methods" in which the matrix $A$ is not required to be strictly lower-triangular. In this situation, the intermediate stages $\{\mathbf{k}_i\}$ exhibit cyclic dependencies and cannot be directly evaluated in a linear fashion. Instead, a nonlinear root-finding algorithm (e.g., Newton's method or fixed-point iteration) is required to find a self-consistent set of values for $\{\mathbf{k}_i\}$. Thus, implicit Runge–Kutta methods are much harder to implement than explicit Runge–Kutta methods. However, implicit methods can achieve much higher order than explicit methods having the same number of stages, and they generally exhibit better stability properties on stiff initial value problems.

methods, respectively, are shown below.

$$
\begin{array}{c|c}
0 & \\
\hline
 & 1
\end{array}
\qquad
\begin{array}{c|cc}
0 & & \\
1/2 & 1/2 & \\
\hline
 & 0 & 1
\end{array}
\qquad
\begin{array}{c|cc}
0 & & \\
1 & 1 & \\
\hline
 & 1/2 & 1/2
\end{array}
\qquad
\begin{array}{c|cccc}
0 & & & & \\
1/2 & 1/2 & & & \\
1/2 & 0 & 1/2 & & \\
1 & 0 & 0 & 1 & \\
\hline
 & 1/6 & 1/3 & 1/3 & 1/6
\end{array}
\tag{1.19}
$$

Note that it is conventional to omit all zeroes on and above the main diagonal of $A$, since these parameters are fixed by definition and convey no meaningful information.

**Remark 1.4.** In this paper, the phrase "Runge–Kutta method" will always be assumed to mean "explicit Runge–Kutta method" unless otherwise specified. (See Footnote 2 for details.)

**Remark 1.5.** Because initial value problems can always be rewritten in autonomous form, the vector $\mathbf{c}$ is not strictly necessary to describe a Runge–Kutta method. By the procedure described in Remark 1.2, knowledge of $A$ and $\mathbf{b}$ is sufficient to apply a Runge–Kutta method to any initial value problem, which effectively sets

$$c_i := \sum_{j=1}^{i-1} A_{ij}. \tag{1.20}$$

Thus, in this paper, we will always treat Runge–Kutta methods as being defined by $(A, \mathbf{b})$, with the value of $\mathbf{c}$ prescribed by Equation (1.20). Under this convention, an $s$-stage Runge–Kutta method is defined by $s(s+1)/2$ real-valued parameters.

Like all numerical approximation algorithms, Runge–Kutta methods are characterized by the amount of computational work they perform and the accuracy of the result they obtain. These properties are captured in two fundamental parameters: *stages* and *order*.

**Definition 1.6.** A ***stage*** in a Runge–Kutta method is a single evaluation of $\mathbf{f}$. For $s \in \mathbb{N}$, an $s$-stage Runge–Kutta method is one that performs $s$ evaluations of $\mathbf{f}$.

**Definition 1.7.** A Runge–Kutta method is said to have ***order*** $p \in \mathbb{N}$ if for every initial value problem $(\mathbf{f}, t_0, \mathbf{y}_0)$ with exact solution $\mathbf{y}$, the result $\mathbf{y}_1$ produced by the method satisfies

$$\|\mathbf{y}_1 - \mathbf{y}(t_0 + h)\| = O(h^{p+1}) \qquad \text{as} \qquad h \to 0. \tag{1.21}$$

The usual goal when designing Runge–Kutta methods is to produce the most accurate result possible using the minimum amount of computational work. In these terms, we aim to design methods that achieve high order using a small number of stages. All four of the methods presented so far are known to be optimal in this sense. Euler's method is the unique example of a 1-stage Runge–Kutta method of order 1, since it follows directly from Taylor expansion of the exact solution that

$$\mathbf{y}(t_0 + h) = \mathbf{y}(t_0) + h\mathbf{y}'(t_0) + O(h^2) = \mathbf{y}_0 + h\,\mathbf{f}(t_0, \mathbf{y}_0) + O(h^2). \tag{1.22}$$

Similar arguments show that the midpoint method and Heun's method are 2-stage methods of order 2, and through laborious computations, that the classical Runge–Kutta method is a 4-stage method of order 4. In fact, the minimum number of stages required is known for all orders $p \leq 8$. (These results are further discussed in Section 1.4.)

| Order | Minimum number of stages |
|---|---|
| 1 | 1 |
| 2 | 2 |
| 3 | 3 |
| 4 | 4 |
| 5 | 6 |
| 6 | 7 |
| 7 | 9 |
| 8 | 11 |

Unfortunately, it is unknown where the optimal balance occurs in general.

**Open Problem 1.8.** For $p > 8$, what is the minimum number of stages required in an explicit Runge–Kutta method of order $p$? Equivalently, for $s > 11$, what is the maximum order achievable by an $s$-stage explicit Runge–Kutta method?

## 1.3 The Numerical Analysis of Initial Value Problems

In order to construct high-order Runge–Kutta methods, we first need a reasonable criterion that determines which Runge–Kutta methods achieve a particular order. The standard approach to determining the order of a Runge–Kutta method is to Taylor-expand the approximate solution $\mathbf{y}_1$ and the exact solution $\mathbf{y}(t_0 + h)$ about $t_0$, then verify that all coefficients on terms of order[3] $h, h^2, \dots, h^p$ coincide. Thus, we will begin our analysis by determining what terms occur in the Taylor expansion of

$$\mathbf{y}(t_0 + h) = \mathbf{y}(t_0) + \mathbf{y}'(t_0)h + \mathbf{y}''(t_0)\frac{h^2}{2!} + \mathbf{y}'''(t_0)\frac{h^3}{3!} + \cdots. \tag{1.23}$$

Observe that an ordinary differential equation of the form

$$\mathbf{y}'(t) = \mathbf{f}(t, \mathbf{y}(t)) \tag{1.24}$$

can be used to express all higher derivatives of $\mathbf{y}$ in terms of derivatives of $\mathbf{f}$. For example, the second derivative of $\mathbf{y}$ can be computed as follows:

$$\begin{aligned} \mathbf{y}''(t) &= \frac{\mathrm{d}}{\mathrm{d}t}\mathbf{y}'(t) \\ &= \frac{\mathrm{d}}{\mathrm{d}t}\mathbf{f}(t, \mathbf{y}(t)) \\ &= \partial_t \mathbf{f}(t, \mathbf{y}(t)) + \big[\nabla_{\mathbf{y}} \mathbf{f}(t, \mathbf{y}(t))\big]\mathbf{y}'(t) \\ &= \partial_t \mathbf{f}(t, \mathbf{y}(t)) + \big[\nabla_{\mathbf{y}} \mathbf{f}(t, \mathbf{y}(t))\big]\mathbf{f}(t, \mathbf{y}(t)) \end{aligned} \tag{1.25}$$

Here, $\big[\nabla_{\mathbf{y}} \mathbf{f}(t, \mathbf{y}(t))\big]$ denotes the Jacobian matrix of $\mathbf{f}$ with respect to $\mathbf{y}$, so the expression $\big[\nabla_{\mathbf{y}} \mathbf{f}(t, \mathbf{y}(t))\big]\mathbf{f}(t, \mathbf{y}(t))$ is a matrix-vector product. Similarly, the third derivative of $\mathbf{y}$ can be computed as follows:

$$\begin{aligned} \mathbf{y}'''(t) &= \frac{\mathrm{d}}{\mathrm{d}t}\mathbf{y}''(t) \\ &= \frac{\mathrm{d}}{\mathrm{d}t}\Big[\partial_t \mathbf{f}(t, \mathbf{y}(t)) + \big[\nabla_{\mathbf{y}} \mathbf{f}(t, \mathbf{y}(t))\big]\mathbf{f}(t, \mathbf{y}(t))\Big] \\ &= \frac{\mathrm{d}}{\mathrm{d}t}[\partial_t \mathbf{f}(t, \mathbf{y}(t))] + \Big[\frac{\mathrm{d}}{\mathrm{d}t}\nabla_{\mathbf{y}} \mathbf{f}(t, \mathbf{y}(t))\Big]\mathbf{f}(t, \mathbf{y}(t)) + \nabla_{\mathbf{y}} \mathbf{f}(t, \mathbf{y}(t))\Big[\frac{\mathrm{d}}{\mathrm{d}t}\mathbf{y}'(t)\Big] \\ &= \partial_t^2 \mathbf{f}(t, \mathbf{y}(t)) + \big[\partial_t \nabla_{\mathbf{y}} \mathbf{f}(t, \mathbf{y}(t))\big]\mathbf{y}'(t) \\ &\quad + \big[\partial_t \nabla_{\mathbf{y}} \mathbf{f}(t, \mathbf{y}(t)) + \big[\nabla_{\mathbf{y}}^2 \mathbf{f}(t, \mathbf{y}(t))\big]\mathbf{y}'(t)\big]\mathbf{f}(t, \mathbf{y}(t)) \\ &\quad + \big[\nabla_{\mathbf{y}} \mathbf{f}(t, \mathbf{y}(t))\big]\big[\partial_t \mathbf{f}(t, \mathbf{y}(t)) + \big[\nabla_{\mathbf{y}} \mathbf{f}(t, \mathbf{y}(t))\big]\mathbf{f}(t, \mathbf{y}(t))\big] \\ &= \partial_t^2 \mathbf{f}(t, \mathbf{y}(t)) + 2\big[\partial_t \nabla_{\mathbf{y}} \mathbf{f}(t, \mathbf{y}(t))\big]\mathbf{f}(t, \mathbf{y}(t)) + \big[\big[\nabla_{\mathbf{y}}^2 \mathbf{f}(t, \mathbf{y}(t))\big]\mathbf{f}(t, \mathbf{y}(t))\big]\mathbf{f}(t, \mathbf{y}(t)) \\ &\quad + \big[\nabla_{\mathbf{y}} \mathbf{f}(t, \mathbf{y}(t))\big]\big[\partial_t \mathbf{f}(t, \mathbf{y}(t))\big] + \big[\nabla_{\mathbf{y}} \mathbf{f}(t, \mathbf{y}(t))\big]\big[\nabla_{\mathbf{y}} \mathbf{f}(t, \mathbf{y}(t))\big]\mathbf{f}(t, \mathbf{y}(t)) \end{aligned} \tag{1.26}$$

Here, $\big[\nabla_{\mathbf{y}}^2 \mathbf{f}(t, \mathbf{y}(t))\big]$ denotes the Hessian of $\mathbf{f}$ with respect to $\mathbf{y}$, a symmetric bilinear map $\mathbb{R}^n \times \mathbb{R}^n \to \mathbb{R}^n$ which is fed two copies of $\mathbf{f}(t, \mathbf{y}(t))$ to produce a vector. These examples illustrate that the complexity of these formulas grows very quickly, since the product and chain rules introduce additional terms at an exponential rate.

[3]Note that no verification is necessary for terms of order $h^0$, since by definition, $\mathbf{y}_1 = \mathbf{y}_0 + O(h)$.

One element of this complexity can be removed if we assume, without loss of generality, that our differential equation is written in *autonomous form*:

$$\mathbf{y}'(t) = \mathbf{f}(\mathbf{y}(t)) \tag{1.27}$$

This allows us to eliminate all $t$-derivatives of $\mathbf{f}$ from consideration, killing all terms that contain $\partial_t$. We will henceforth denote derivatives of $\mathbf{f}$ with respect to $\mathbf{y}$ using primes, as no confusion can occur with $t$ eliminated. For notational simplicity, we will also suppress the arguments of $\mathbf{f}$ and its derivatives, with the understanding that every occurrence of $\mathbf{f} = \mathbf{f}(t, \mathbf{y}(t))$ is to be evaluated at the same point. Thus, we will write:

$$\big[\nabla_{\mathbf{y}} \mathbf{f}(t, \mathbf{y}(t))\big] \mathbf{f}(t, \mathbf{y}(t)) \longrightarrow \mathbf{f}'(\mathbf{f}) \tag{1.28}$$

$$\big[\big[\nabla_{\mathbf{y}}^2 \mathbf{f}(t, \mathbf{y}(t))\big] \mathbf{f}(t, \mathbf{y}(t))\big] \mathbf{f}(t, \mathbf{y}(t)) \longrightarrow \mathbf{f}''(\mathbf{f}, \mathbf{f}) \tag{1.29}$$

$$\big[\nabla_{\mathbf{y}} \mathbf{f}(t, \mathbf{y}(t))\big]\big[\nabla_{\mathbf{y}} \mathbf{f}(t, \mathbf{y}(t))\big] \mathbf{f}(t, \mathbf{y}(t)) \longrightarrow \mathbf{f}'(\mathbf{f}'(\mathbf{f})) \tag{1.30}$$

Adopting these notational conventions allows us to rewrite Equations (1.24), (1.25), and (1.26) in a much more palatable form:

$$\mathbf{y}' = \mathbf{f} \tag{1.31}$$

$$\mathbf{y}'' = \mathbf{f}'(\mathbf{f}) \tag{1.32}$$

$$\mathbf{y}''' = \mathbf{f}''(\mathbf{f}, \mathbf{f}) + \mathbf{f}'(\mathbf{f}'(\mathbf{f})) \tag{1.33}$$

In addition, we can write down reasonable formulas for higher derivatives that would have been hideously complicated in the old notation.

$$\mathbf{y}'''' = \mathbf{f}'''(\mathbf{f}, \mathbf{f}, \mathbf{f}) + 3\mathbf{f}''(\mathbf{f}'(\mathbf{f}), \mathbf{f}) + \mathbf{f}'(\mathbf{f}''(\mathbf{f}, \mathbf{f})) + \mathbf{f}'(\mathbf{f}'(\mathbf{f}'(\mathbf{f}))) \tag{1.34}$$

$$\begin{aligned} \mathbf{y}''''' = {} & \mathbf{f}''''(\mathbf{f}, \mathbf{f}, \mathbf{f}, \mathbf{f}) + 6\mathbf{f}'''(\mathbf{f}'(\mathbf{f}), \mathbf{f}, \mathbf{f}) + 4\mathbf{f}''(\mathbf{f}''(\mathbf{f}, \mathbf{f}), \mathbf{f}) \\ & + 4\mathbf{f}''(\mathbf{f}'(\mathbf{f}'(\mathbf{f})), \mathbf{f}) + 3\mathbf{f}''(\mathbf{f}'(\mathbf{f}), \mathbf{f}'(\mathbf{f})) + \mathbf{f}'(\mathbf{f}'''(\mathbf{f}, \mathbf{f}, \mathbf{f})) \\ & + 3\mathbf{f}'(\mathbf{f}''(\mathbf{f}'(\mathbf{f}), \mathbf{f})) + \mathbf{f}'(\mathbf{f}'(\mathbf{f}''(\mathbf{f}, \mathbf{f}))) + \mathbf{f}'(\mathbf{f}'(\mathbf{f}'(\mathbf{f}'(\mathbf{f})))) \end{aligned} \tag{1.35}$$

Here, the $k^{\text{th}}$ derivative of $\mathbf{f}$ is interpreted as a symmetric multilinear function

$$\mathbf{f}^{(k)} : \underbrace{\mathbb{R}^n \times \cdots \times \mathbb{R}^n}_{k \text{ copies}} \to \mathbb{R}^n. \tag{1.36}$$

(Symmetry of $\mathbf{f}^{(k)}$ is a consequence of the fact that partial derivative operators commute.) Notice that some of these terms appear with an integer coefficient greater than one because they arise from lower-order terms in more than one way. For example, $\mathbf{f}''(\mathbf{f}'(\mathbf{f}), \mathbf{f})$ is obtained three times: twice from differentiating $\mathbf{f}''(\mathbf{f}, \mathbf{f})$ with respect to each inner copy of $\mathbf{f}$, and once from differentiating $\mathbf{f}'(\mathbf{f}'(\mathbf{f}))$ with respect to the outer copy of $\mathbf{f}'$. Symmetry allows us to identify $\mathbf{f}''(\mathbf{f}'(\mathbf{f}), \mathbf{f}) = \mathbf{f}''(\mathbf{f}, \mathbf{f}'(\mathbf{f}))$.

The functions $\mathbb{R}^n \to \mathbb{R}^n$ that arise in this manner, consisting of compositions of $\mathbf{f}$ with its derivatives, are called the ***elementary differentials*** of $\mathbf{f}$. Cayley first observed, over 160 years ago [18], that the structures of elementary differentials are precisely captured by (isomorphism classes of) rooted trees. These notions are formally defined below.

**Definition 1.9.** A ***rooted tree*** is an ordered pair $\tau = (G, v)$ consisting of a connected acyclic undirected finite graph $G$ and a distinguished vertex $v \in V(G)$, called the ***root***. As with graphs, the number of vertices in a rooted tree $\tau$ is called its ***order*** and denoted by $|\tau|$. Two rooted trees are said to be ***isomorphic*** if there exists an isomorphism of their underlying graphs that sends the root of one tree to the root of the other. For $n \in \mathbb{N}$, we denote by $T_n$ the set of all isomorphism classes of rooted trees of order $n$.

We will henceforth identify the phrases "rooted tree" and "isomorphism class of rooted trees," since for our purposes, this causes no confusion. Rooted trees will always be drawn with their roots at the bottom. To illustrate this convention, all rooted trees of order $\leq 4$ are depicted below.

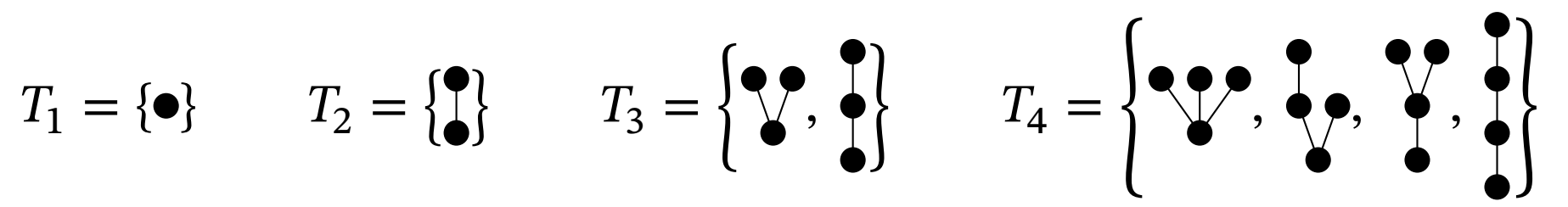

(1.37)

To formalize the correspondence between rooted trees and elementary differentials, we will first need a notion of "building up" a rooted tree from its subtrees.

**Definition 1.10.** Let $\tau_1, \tau_2, \ldots, \tau_n$ be rooted trees and $k_1, k_2, \ldots, k_n \in \mathbb{N}$. We denote by $[\tau_1^{k_1}, \tau_2^{k_2}, \ldots, \tau_n^{k_n}]$ the rooted tree obtained by taking the disjoint union of $k_1$ copies of $\tau_1$, $k_2$ copies of $\tau_2$, ..., $k_n$ copies of $\tau_n$, and adjoining a new vertex which is adjacent to each root in each copy of $\tau_1, \tau_2, \ldots, \tau_n$. This new vertex is then declared the root of $[\tau_1^{k_1}, \tau_2^{k_2}, \ldots, \tau_n^{k_n}]$.

Several examples are presented below to illustrate this definition.

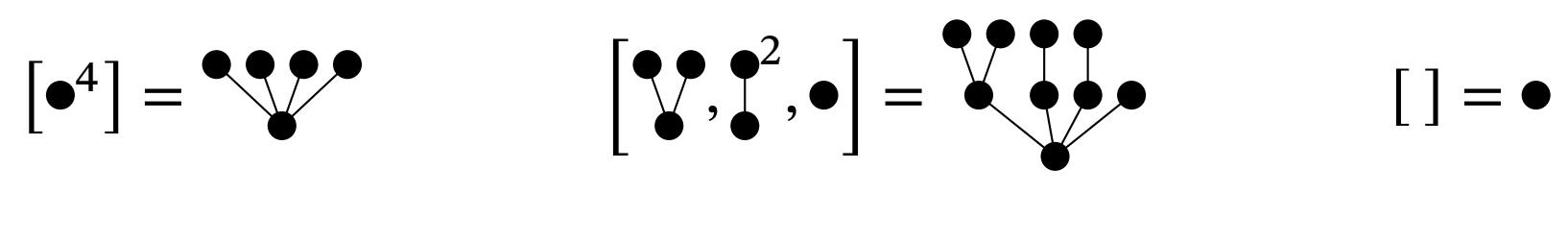

(1.38)

Note that in the final example, we take the disjoint union of an empty sequence of rooted trees, obtaining a graph with zero vertices, then adjoin a new vertex to obtain the unique one-vertex rooted tree. In fact, the observation that $[\,] = \bullet$ allows us to write any rooted tree purely in terms of brackets.

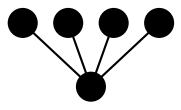 $= [[\,],[\,],[\,],[\,]]$ 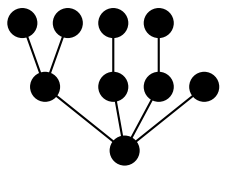 $= [[[\,],[\,]],[[\,]],[[\,]],[\,]]$ (1.39)

In addition, observe that any rooted tree $\tau$ can be written as the bracket of its ***legs***, which are the rooted subtrees obtained from $\tau$ by deleting its root and taking the connected components of the resulting forest. The root in each component is the unique vertex which was previously adjacent to the root of $\tau$. This allows us to write any rooted tree as $\tau = [\tau_1, \tau_2, \dots, \tau_n]$ without loss of generality. Moreover, by partitioning the legs into isomorphism classes, we can also write $\tau = [\tau_1^{k_1}, \tau_2^{k_2}, \dots, \tau_n^{k_n}]$ where $\tau_1, \tau_2, \dots, \tau_n$ are pairwise non-isomorphic. Both decompositions will be useful in stating recursive definitions of functions on rooted trees. Note that recursive definitions of this type are always well-founded because legs always have strictly fewer vertices than the rooted tree they came from.

**Definition 1.11.** Let $\mathbf{f} : \mathbb{R}^n \to \mathbb{R}^n$ be a smooth function, and let $\tau$ be a rooted tree. The ***elementary differential*** of $\mathbf{f}$ corresponding to $\tau$ is the function $D_{\mathbf{f}}(\tau) : \mathbb{R}^n \to \mathbb{R}^n$ defined recursively as follows:

$$D_{\mathbf{f}}(\bullet) := \mathbf{f} \tag{1.40}$$

$$D_{\mathbf{f}}([\tau_1, \tau_2, \dots, \tau_m]) := \mathbf{f}^{(m)}(D_{\mathbf{f}}(\tau_1), D_{\mathbf{f}}(\tau_2), \dots, D_{\mathbf{f}}(\tau_m)) \tag{1.41}$$

Several examples are presented below to illustrate this definition.

$D_{\mathbf{f}}\Big($ 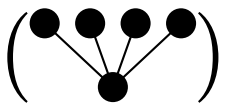 $\Big) = \mathbf{f}''''(\mathbf{f}, \mathbf{f}, \mathbf{f}, \mathbf{f})$ $D_{\mathbf{f}}\Big($ 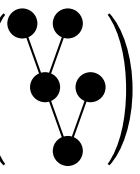 $\Big) = \mathbf{f}''(\mathbf{f}''(\mathbf{f}, \mathbf{f}), \mathbf{f})$

$D_{\mathbf{f}}\Big($ 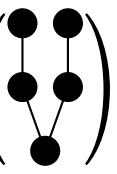 $\Big) = \mathbf{f}''(\mathbf{f}'(\mathbf{f}), \mathbf{f}'(\mathbf{f}))$ $D_{\mathbf{f}}\Big($ 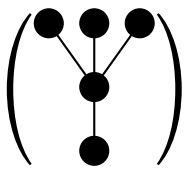 $\Big) = \mathbf{f}'(\mathbf{f}'''(\mathbf{f}, \mathbf{f}, \mathbf{f}))$ (1.42)

For further illustration, the derivative formulas (1.31)–(1.35) are rewritten below using rooted trees.

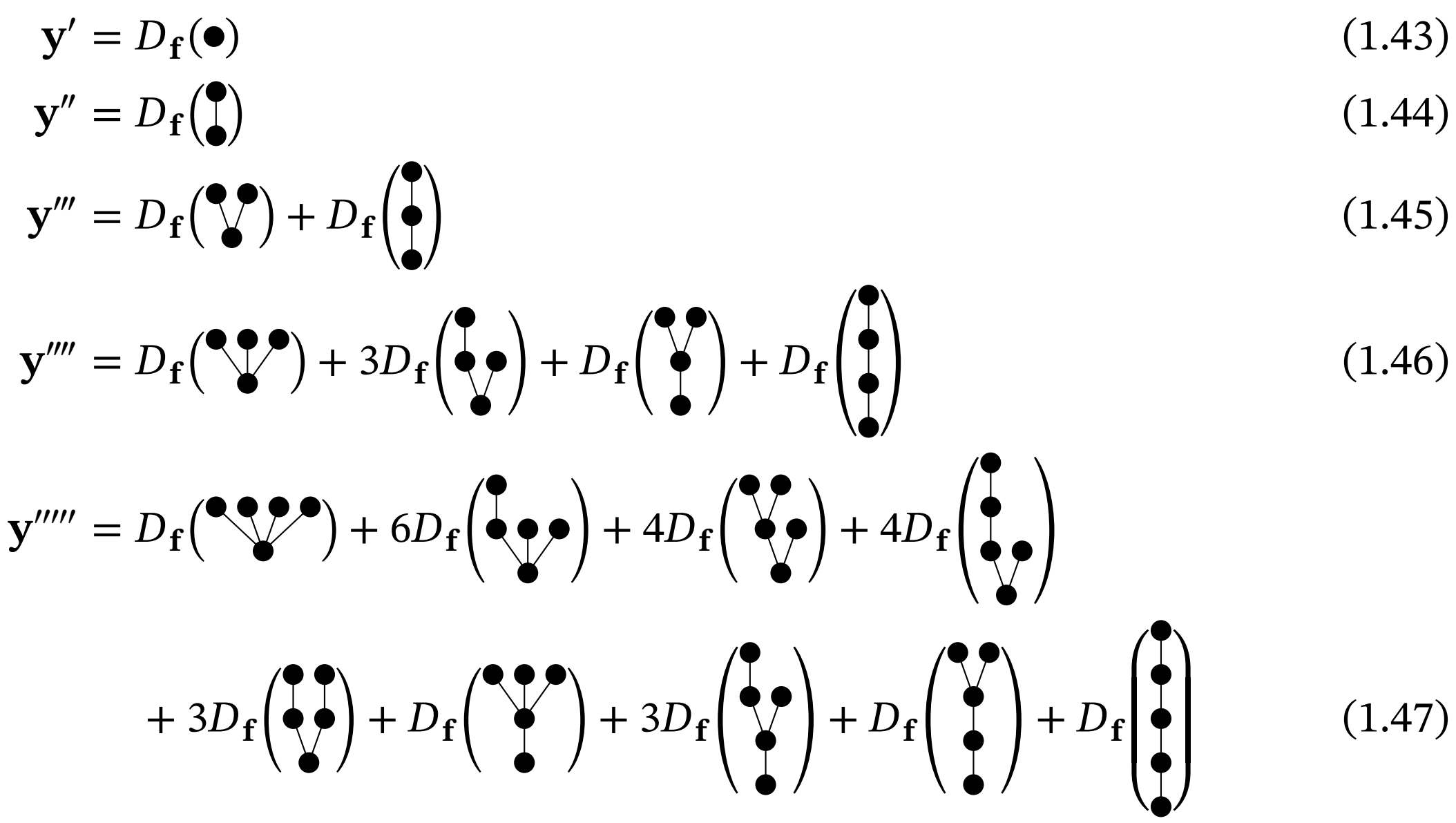

The coefficients that appear in front of each elementary differential in these formulas turn out to be combinatorial quantities that can be directly computed from the corresponding rooted tree $\tau$.

**Definition 1.12.** Let $\tau$ be a rooted tree, and let $D$ be the directed graph obtained from $\tau$ by orienting all edges away from the root. A ***rooted labeling*** of $\tau$ is a bijection $\ell : V(D) \to \{1, 2, \dots, |\tau|\}$ having the following property: for all $v, w \in V(D)$, if there exists a directed edge $v \to w$, then $\ell(v) < \ell(w)$. Two rooted labelings of $\tau$ are said to be ***equivalent*** if one can be written as the composition of the other with an automorphism of $\tau$. We denote by $\alpha(\tau)$ the number of equivalence classes of rooted labelings of $\tau$.

**Proposition 1.13.** Let $\mathbf{y} : \mathbb{R} \to \mathbb{R}^n$ and $\mathbf{f} : \mathbb{R}^n \to \mathbb{R}^n$ be smooth functions satisfying $\mathbf{y}'(\tau) = \mathbf{f}(\mathbf{y}(\tau))$ for all $t \in \mathbb{R}$. For all $k \in \mathbb{N}$, the $k^{\text{th}}$ derivative of $\mathbf{y}$ is given by

$$\mathbf{y}^{(k)} = \sum_{\tau \in T_k} \alpha(\tau) D_{\mathbf{f}}(\tau). \tag{1.48}$$

*Proof sketch.* For $k = 1$, the claim $\mathbf{y}' = \alpha(\bullet) D_{\mathbf{f}}(\bullet) = \mathbf{f}$ holds by hypothesis, so suppose $k > 1$. A rooted labeling of a rooted tree $\tau \in T_k$ can be thought of a specifying a sequence of instructions for building $\tau$ by starting from a single vertex and repeatedly attaching leaves. The number of non-equivalent ways to do this coincides with the number of times $D_{\mathbf{f}}(\tau)$ appears when starting from the equation $\mathbf{y}' = \mathbf{f}$ and differentiating both sides $k - 1$ times. $\square$

With Proposition 1.13 in hand, we can now state the correct analogue of Taylor's Theorem for initial value problems.

**Corollary 1.14.** Let $(\mathbf{f}, t_0, \mathbf{y}_0)$ be an initial value problem with solution $\mathbf{y}$. For all $p \in \mathbb{N}$,

$$\mathbf{y}(t_0 + h) = \mathbf{y}_0 + \sum_{k=1}^{p} \sum_{\tau \in T_k} \frac{\alpha(\tau) h^k}{k!} D_{\mathbf{f}}(\tau)(\mathbf{y}_0) + O(h^{p+1}) \qquad \text{as} \qquad h \to 0. \tag{1.49}$$

*Proof.* Plug the rooted tree expansion (1.48) into the usual statement of Taylor's theorem. □

Unfortunately, this analogue of Taylor expansion turns out to be cumbersome, since the number of rooted trees grows exponentially in the number of vertices. The precise asymptotic behavior is stated below, confirming our earlier suspicion that "the product and chain rules introduce additional terms at an exponential rate."

**Proposition 1.15** (Otter 1948 [48])**.** The number of isomorphism classes of rooted trees of order $n$ satisfies

$$|T_n| = \Theta\left(\frac{\alpha^n}{n^{3/2}}\right) \qquad \text{as} \qquad n \to \infty, \tag{1.50}$$

where $\alpha = 2.955765285652\ldots$ denotes the ***rooted tree constant***.

To illustrate this rapid growth, a table of values of $|T_n|$ and the cumulative sums $|T_1|+|T_2|+\cdots+|T_n|$ is shown on the following page. By Corollary 1.14, $|T_1|+|T_2|+\cdots+|T_n|$ is the number of terms in the $n^{\text{th}}$-order Taylor expansion of $\mathbf{y}(t_0 + h)$.

Given this exponential rate of growth, we might naturally wonder if it is really necessary to separately consider the elementary differentials corresponding to every rooted tree. After all, it is conceivable that Equation (1.49) could be simplified by writing some elementary differential as a linear combination of others. Unfortunately, the following result demonstrates that this is impossible; in general, all elementary differentials are linearly independent.

**Proposition 1.16** (Butcher 2003 [16], pp. 146–147)**.** Let $T = \{\tau_1, \tau_2, \ldots, \tau_n\}$ be a finite set of rooted trees. There exists an $n$-dimensional autonomous initial value problem $(\mathbf{f}, t_0, \mathbf{y}_0)$ in which the elementary differentials of $\mathbf{f}$ satisfy

$$D_{\mathbf{f}}(\tau_i)(\mathbf{y}_0) = \mathbf{e}_i \qquad \text{for all} \qquad i = 1, 2, \ldots, n, \tag{1.51}$$

where $\mathbf{e}_i$ denotes the $i^{\text{th}}$ stzandard basis vector in $\mathbb{R}^n$.

| $n$ | $\lvert T_n \rvert$ | $\lvert T_1 \rvert + \lvert T_2 \rvert + \cdots + \lvert T_n \rvert$ |
|---|---|---|
| 1 | 1 | 1 |
| 2 | 1 | 2 |
| 3 | 2 | 4 |
| 4 | 4 | 8 |
| 5 | 9 | 17 |
| 6 | 20 | 37 |
| 7 | 48 | 85 |
| 8 | 115 | 200 |
| 9 | 286 | 486 |
| 10 | 719 | 1205 |
| 11 | 1842 | 3047 |
| 12 | 4766 | 7813 |

## 1.4 Order Conditions and Butcher Weights

Having determined the Taylor expansion of the exact solution $\mathbf{y}(t_0 + h)$ about $t_0$, it now remains to compute the Taylor expansion of the approximate solution $\mathbf{y}_1$ produced by a Runge–Kutta method defined by arbitrary parameters $(A, \mathbf{b})$. The basic idea is to inductively compute the Taylor expansion of the intermediate stages (written here in autonomous form)

$$\mathbf{k}_i := h\,\mathbf{f}\left(\mathbf{y}_0 + \sum_{j=1}^{i-1} a_{ij}\mathbf{k}_j\right) \qquad \text{for} \qquad i = 1, 2, \dots, s \tag{1.52}$$

about $\mathbf{y}_0$, then take dot products with $\mathbf{b}$ to obtain the Taylor expansion of

$$\mathbf{y}_1 := \mathbf{y}_0 + \sum_{i=1}^{s} b_i \mathbf{k}_i. \tag{1.53}$$

Beyond this initial idea, the majority of the argument consists of tracking indices across invocations of the inductive hypothesis. This bookkeeping is relatively unenlightening, so we will develop only the necessary machinery to state the main result, and refer the interested reader to Section 313 of Butcher's textbook [16] for details of the proof.

**Definition 1.17.** Let $\tau$ be a rooted tree and $A \in \mathbb{R}^{s \times s}$. The ***Butcher weight*** of $A$

corresponding to $\tau$ is the vector $\mathbf{\Phi}_A(\tau) \in \mathbb{R}^s$ recursively defined by

$$\mathbf{\Phi}_A(\bullet) := \mathbf{1} \tag{1.54}$$

$$\mathbf{\Phi}_A([\tau_1, \tau_2, \ldots, \tau_m]) := (A\mathbf{\Phi}_A(\tau_1)) \odot (A\mathbf{\Phi}_A(\tau_2)) \odot \cdots \odot (A\mathbf{\Phi}_A(\tau_m)) \tag{1.55}$$

where $\mathbf{1} \in \mathbb{R}^s$ denotes the $s$-dimensional vector whose entries are all 1, and $\odot$ denotes the elementwise product of vectors.

Several examples are presented below to illustrate this definition. In accordance with the convention established in Remark 1.5, we will henceforth write $\mathbf{c} := A\mathbf{1}$.

$\mathbf{\Phi}_A\Big($ 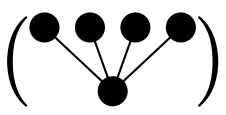 $\Big) = \mathbf{c}^{\odot 4}$ $\qquad$ $\mathbf{\Phi}_A\Big($ 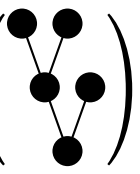 $\Big) = (A(\mathbf{c}^{\odot 2})) \odot (A\mathbf{c})$

$\mathbf{\Phi}_A\Big($ 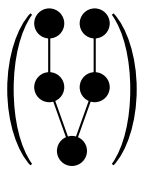 $\Big) = (A\mathbf{c})^{\odot 2}$ $\qquad$ $\mathbf{\Phi}_A\Big($ 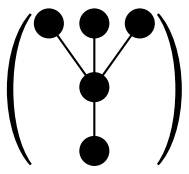 $\Big) = A(\mathbf{c}^{\odot 3})$ $\qquad$ (1.56)

Here, $\mathbf{v}^{\odot n}$ denotes the elementwise $n^{\text{th}}$ power of $\mathbf{v}$. For example, $\mathbf{v}^{\odot 3} := \mathbf{v} \odot \mathbf{v} \odot \mathbf{v}$.

**Definition 1.18.** The ***symmetry*** of a rooted tree $\tau$, denoted by $\sigma(\tau)$, is the order of the automorphism group of $\tau$.

The following proposition gives a straightforward recursive algorithm for computing the symmetry of a rooted tree.

**Proposition 1.19.** Let $\tau = [\tau_1^{k_1}, \tau_2^{k_2}, \ldots, \tau_n^{k_n}]$ be a rooted tree with distinct legs $\tau_1, \tau_2, \ldots, \tau_n$. Then

$$\sigma(\tau) = \prod_{i=1}^{n} k_i!\sigma(\tau_i). \tag{1.57}$$

*Proof.* Observe that each automorphism of $\tau$ can be written uniquely as a permutation of its legs followed by an automorphism within each leg. □

**Lemma 1.20** (Butcher 1963 [11]). Let $(\mathbf{f}, t_0, \mathbf{y}_0)$ be an initial value problem, $A \in \mathbb{R}^{s \times s}$ a strictly lower-triangular matrix, $\mathbf{b} \in \mathbb{R}^s$, and $h \in \mathbb{R}$. Let $\mathbf{y}_1$ denote the result of applying the Runge–Kutta method defined by $(A, \mathbf{b})$ to $(\mathbf{f}, t_0, \mathbf{y}_0)$ with step size $h$. For all $p \in \mathbb{N}$,

$$\mathbf{y}_1 = \mathbf{y}_0 + \sum_{k=1}^{p} \sum_{\tau \in T_k} \frac{(\mathbf{b} \cdot \mathbf{\Phi}_A(\tau))h^k}{\sigma(\tau)} D_{\mathbf{f}}(\tau)(\mathbf{y}_0) + O(h^{p+1}) \quad \text{as } h \to 0. \tag{1.58}$$

With Corollary 1.14 and Lemma 1.20 in hand, we can now compare the Taylor expansions of the exact (1.49) and approximate (1.58) solutions term-by-term. We see that for each rooted tree $\tau$, the coefficient of $D_{\mathbf{f}}(\tau)(\mathbf{y}_0)$ in the exact expansion is $\alpha(\tau)h^{|\tau|}/|\tau|!$, while the corresponding coefficient in the approximate expansion is $(\mathbf{b} \cdot \mathbf{\Phi}_A(\tau))h^{|\tau|}/\sigma(\tau)$. This comparison naturally leads us to seek a relationship between the quantities $|\tau|!$, $\alpha(\tau)$, and $\sigma(\tau)$. Luckily, such a relationship exists through a fourth quantity introduced in the following definition.

**Definition 1.21.** Let $\tau = [\tau_1, \tau_2, \dots, \tau_n]$ be a rooted tree. The ***factorial*** or ***density*** of $\tau$ is the number $\tau! \in \mathbb{N}$ recursively defined by:

$$(\bullet)! := 1 \tag{1.59}$$

$$[\tau_1, \tau_2, \dots, \tau_n]! := \left(1 + \sum_{i=1}^{n} |\tau_i|\right) \prod_{i=1}^{n} \tau_i! \tag{1.60}$$

This recursive definition can be understood as labeling each vertex of $\tau$ with an integer corresponding to the order of the subtree of $\tau$ rooted at that vertex. The factorial of $\tau$ is the product of all such labels, as illustrated by the following example.

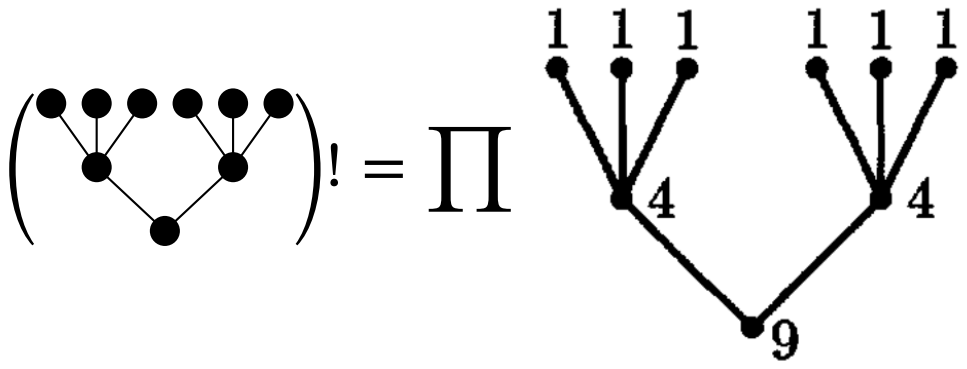

$$= 9 \cdot 4^2 \cdot 1^6 = 144 \tag{1.61}$$

**Proposition 1.22.** For all rooted trees $\tau$,

$$\alpha(\tau) = \frac{|\tau|!}{\tau!\sigma(\tau)}. \tag{1.62}$$

*Proof sketch.* Observe that $|\tau|!$ is the number of unrestricted labelings of $\tau$ (i.e., bijections $V(\tau) \to \{1, 2, \dots, |\tau|\}$). Restricting attention to rooted labelings of $\tau$ precisely means requiring that the label applied to each vertex $v \in V(\tau)$ be the smallest label in the subtree of $\tau$ rooted at $v$. This reduces the number of labelings by a factor of $\tau!$. Each equivalence class of rooted labelings has $\sigma(\tau)$ members, so the number of equivalence classes of rooted labelings is precisely $\alpha(\tau) = |\tau|!/(\tau!\sigma(\tau))$. □

This relationship between $\alpha(\tau)$, $\sigma(\tau)$, $|\tau|!$, and $\tau!$ allows us to state the order conditions for Runge–Kutta methods in a particularly elegant form.

**Theorem 1.23** (Runge–Kutta Order Conditions)**.** Let $A \in \mathbb{R}^{s \times s}$ and $\mathbf{b} \in \mathbb{R}^s$. The

Runge–Kutta method defined by $(A, \mathbf{b})$ has order $p \in \mathbb{N}$ if and only if

$$\mathbf{b} \cdot \mathbf{\Phi}_A(\tau) = \frac{1}{\tau!} \qquad \text{for all rooted trees } \tau \text{ of order} \leq p. \tag{1.63}$$

*Proof.* By comparing the Taylor expansions (1.49), (1.58) of the exact solution $\mathbf{y}(t_0 + h)$ and the approximate solution $\mathbf{y}_1$, we see that $\|\mathbf{y}_1 - \mathbf{y}(t_0 + h)\| = O(h^{p+1})$ as $h \to 0$ if and only if the coefficients of $D_\mathbf{f}(\tau)(\mathbf{y}_0)$ agree for all rooted trees $\tau$ of order $\leq p$. (The "only if" implication in the preceding statement follows from the linear independence of elementary differentials established in Proposition 1.16.) Thus, we require:

$$\frac{(\mathbf{b} \cdot \mathbf{\Phi}_A(\tau))h^{|\tau|}}{\sigma(\tau)} = \frac{\alpha(\tau)h^{|\tau|}}{|\tau|!}$$

$$\mathbf{b} \cdot \mathbf{\Phi}_A(\tau) = \frac{\alpha(\tau)\sigma(\tau)}{|\tau|!} = \frac{1}{\tau!} \qquad \square$$

To illustrate the typical form of Runge–Kutta order conditions, the equation $\mathbf{b} \cdot \mathbf{\Phi}_A(\tau) = 1/\tau!$ is written explicitly below for all rooted trees $\tau$ of order $\leq 4$. Again, in accordance with the convention established in Remark 1.5, we write $\mathbf{c} \coloneqq A\mathbf{1}$.

$(\bullet)$ : $\mathbf{b} \cdot \mathbf{1} = 1$ 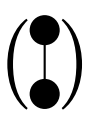 : $\mathbf{b} \cdot \mathbf{c} = \frac{1}{2}$ (1.64)

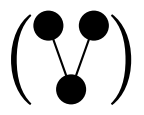 : $\mathbf{b} \cdot (\mathbf{c}^{\odot 2}) = \frac{1}{3}$ 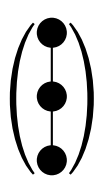 : $\mathbf{b} \cdot (A\mathbf{c}) = \frac{1}{6}$ (1.65)

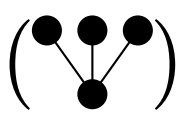 : $\mathbf{b} \cdot (\mathbf{c}^{\odot 3}) = \frac{1}{4}$ 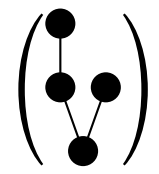 : $\mathbf{b} \cdot ((A\mathbf{c}) \odot \mathbf{c}) = \frac{1}{8}$ (1.66)

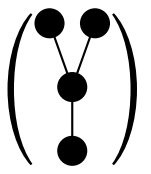 : $\mathbf{b} \cdot (A(\mathbf{c}^{\odot 2})) = \frac{1}{12}$ 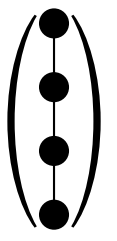 : $\mathbf{b} \cdot (A^2\mathbf{c}) = \frac{1}{24}$ (1.67)

Given the recursive structure of the Butcher weights $\mathbf{\Phi}_A(\tau)$ and the tree factorial $\tau!$, we might naturally wonder if some of the Runge–Kutta order conditions (1.63) are redundant. Unfortunately, the following proposition demonstrates that this is not the case.

**Proposition 1.24** (Butcher 2003 [16], pp. 149–150)**.** Let $T = \{\tau_1, \tau_2, \dots, \tau_n\}$ be a finite set of rooted trees, and let $x_1, x_2, \dots, x_n \in \mathbb{R}$ be arbitrary. There exists an explicit Runge–Kutta method $(A, \mathbf{b})$ for which

$$\mathbf{b} \cdot \mathbf{\Phi}_A(\tau_i) = x_i \qquad \text{for all } i = 1, 2, \dots, n. \tag{1.68}$$

The upshot of this proposition is that we can always cook up a Runge–Kutta method that satisfies any desired set of order conditions (1.63) while violating any other desired set (provided that both sets are finite and disjoint). Thus, the order conditions for a given order $p$ are all necessary; no set of such conditions implies any other, so we cannot eliminate any of them.

Corollary 1.14 and Lemma 1.20 also allow us to estimate the local truncation error committed by a Runge–Kutta method by calculating the leading-order differences between the exact (1.49) and approximate (1.58) Taylor expansions.

**Definition 1.25.** The ***principal error coefficients*** of a Runge–Kutta method $(A, \mathbf{b}, \mathbf{c})$ of order $p$ are the following multiset of real numbers:

$$\left\{ \frac{1}{\sigma(\tau)} \left( \mathbf{b} \cdot \mathbf{\Phi}_A(\tau) - \frac{1}{\tau!} \right) : \tau \in T_{p+1} \right\} \tag{1.69}$$

We can think of the principal error coefficients of a Runge–Kutta method of order $p$ as measuring how badly the method fails to achieve order $p+1$. The normalization factor $1/\sigma(\tau)$ makes the principal error coefficients precisely the coefficients of the leading-order terms $h^{p+1} D_{\mathbf{f}}(\tau)(\mathbf{y}_0)$ in the Taylor expansion of $\mathbf{y}_1 - \mathbf{y}(t_0 + h)$.

## 1.5 Order Barriers and Upper Bounds

To conclude Chapter 1, we demonstrate how the structure of the Runge–Kutta order conditions (1.63) can be used to prove lower bounds on the number of stages required to achieve a given order.

**Theorem 1.26.** For all $p \in \mathbb{N}$, any explicit Runge–Kutta method of order $p$ must have at least $s \geq p$ stages.

*Proof.* We proceed by proving the contrapositive. Suppose $s < p$, and let $\tau = [[\cdots[\,]\cdots]]$ be the rooted tree consisting of a path on $p$ vertices rooted at one of its endpoints. Observe that the order condition associated to $\tau$ is the equation

$$\mathbf{b} \cdot (A^{p-1} \mathbf{1}) = \frac{1}{p!}. \tag{1.70}$$

By definition, $A$ is a strictly lower-triangular $s \times s$ matrix, so $A$ is nilpotent of index $s$. Hence, the assumption $p - 1 \geq s$ implies that $A^{p-1}$ is the zero matrix, and $\mathbf{b} \cdot (A^{p-1}\mathbf{1}) = 0 \neq 1/p!$. □

Unfortunately, this bound is only optimal for $1 \leq p \leq 4$, and extremely complicated arguments are required to prove the known optimal bounds for $5 \leq p \leq 8$. These are the celebrated ***Butcher Barrier Theorems***, proven by John Butcher over the course of several decades, and stated below without proof.

**Theorem 1.27** (Butcher 1964 [12], 1965 [13], 1985 [15]).

- Any explicit Runge–Kutta method of order $p \geq 5$ must have $s \geq p + 1$ stages.
- Any explicit Runge–Kutta method of order $p \geq 7$ must have $s \geq p + 2$ stages.
- Any explicit Runge–Kutta method of order $p \geq 8$ must have $s \geq p + 3$ stages.
- For $5 \leq p \leq 8$, all of these bounds are achieved. That is, there exist a 6-stage method of order 5, a 7-stage method of order 6, a 9-stage method of order 7, and an 11-stage method of order 8.

We warn the interested reader that the proofs of these bounds require intricate casework that meticulously dissects the structure of the Runge–Kutta order conditions to derive algebraic contradictions. It is not clear whether these strategies can be adapted to prove similar bounds for $p \geq 9$, where nothing is currently known. We note that Khashin's discovery of a 13-stage method of order 9 [41] reduces the uncertainty for $p = 9$ to a single stage $12 \leq s \leq 13$, but even for $p = 10$, the gap $13 \leq s \leq 16$ seems far more difficult to cross.

In the opposite direction, there is a straightforward construction that produces Runge–Kutta methods of arbitrary order using only a quadratic number of stages.

**Proposition 1.28.** For all $p \in \mathbb{N}$, there exists an explicit Runge–Kutta method of order $p$ having $s = p(p-1)/2 + 1$ stages.

*Proof sketch.* Apply Richardson extrapolation to $p$ independent runs of Euler's method using step sizes $h, h/2, \dots, h/p$. This process can be written as a Runge–Kutta method having a Butcher table of the following form (illustrated for $p = 4$):

$$\begin{array}{c|ccccccc}
0 & & & & & & & \\
1/2 & 1/2 & & & & & & \\
1/3 & 1/3 & 0 & & & & & \\
2/3 & 1/3 & 0 & 1/3 & & & & \\
1/4 & 1/4 & 0 & 0 & 0 & & & \\
2/4 & 1/4 & 0 & 0 & 0 & 1/4 & & \\
3/4 & 1/4 & 0 & 0 & 0 & 1/4 & 1/4 & \\
\hline
 & 0 & 2 & -9/2 & -9/2 & 8/3 & 8/3 & 8/3
\end{array} \tag{1.71}$$

By construction, this method achieves order $p$ using $s = p(p-1)/2 + 1$ stages. □

We note that extrapolated Runge–Kutta methods of this type typically contain an excessive number of stages compared to methods that are specifically constructed to achieve a particular order. They can also experience catastrophic cancellation due to the coefficients in **b** exhibiting large magnitudes with alternating signs (a characteristic of Richardson extrapolation). Thus, they perform poorly in practice compared to methods which lack these features.

## Chapter 2

# Search Methodology

As shown in Chapter 1, the problem of constructing an $s$-stage Runge–Kutta method of order $p$ reduces to solving a system of $|T_1| + |T_2| + \cdots + |T_p|$ polynomial equations of the form

$$\mathbf{b} \cdot \mathbf{\Phi}_A(\tau) = \frac{1}{\tau!} \qquad \text{for all rooted trees } \tau \text{ of order} \leq p. \tag{2.1}$$

This is a system of equations in $s(s+1)/2$ variables represented by a strictly lower-triangular matrix $A \in \mathbb{R}^{s \times s}$ and a vector $\mathbf{b} \in \mathbb{R}^s$. Proposition 1.15 shows that the number of equations grows exponentially in $p$, while the number of variables grows quadratically in $s$. Thus, for comparable values of $p$ and $s$, the system of Runge–Kutta order conditions quickly becomes strongly overdetermined. Typical situations for optimal ($1 \leq p \leq 8$) and best-known ($9 \leq p \leq 10$) values of $s$ are illustrated in the table below.

| $p$ | $s$ | Number of equations | Number of variables |
|---|---|---|---|
| 1 | 1 | 1 | 1 |
| 2 | 2 | 2 | 3 |
| 3 | 3 | 4 | 6 |
| 4 | 4 | 8 | 10 |
| 5 | 6 | 17 | 21 |
| 6 | 7 | 37 | 28 |
| 7 | 9 | 85 | 45 |
| 8 | 11 | 200 | 66 |
| 9 | 13 | 486 | 91 |
| 10 | 16 | 1205 | 136 |

In some sense, it is a miracle that Runge–Kutta methods of order $p \geq 6$ exist at all. In fact, something even more miraculous is true: all known Runge–Kutta methods of order $p \geq 6$ are not isolated solutions of Equation (2.1), but lie on *positive-dimensional* solution varieties [41]. Thus, they are not only solutions of a strongly overdetermined system of equations, but solutions with free parameters! To the best of our knowledge, no satisfactory explanation of this phenomenon is presently known, but it can be used to our advantage by adjusting the free parameters to reduce the magnitudes of the principal error coefficients.

The primary contribution of this paper is the family of 16-stage Runge–Kutta methods of order 10 presented in Appendix A. These methods lie on a locally 11-dimensional solution variety, of which we present two members: the first such method found by our search procedure (A.1), followed by an optimized method whose principal error coefficients are roughly an order of magnitude smaller (A.2). This chapter describes the numerical techniques by which these methods were discovered. Some of these are well-known, while others are original developments specific to the structure of the Runge–Kutta order conditions. We note that, to the best of our knowledge, this is the first time each of these ideas has been applied to the discovery of Runge–Kutta methods.

All algorithms and techniques described in this chapter have been implemented in an open-source software package called "RKTK: A Runge–Kutta Toolkit." The source code for RKTK, written in a mixture of the C++ [58] and Julia [6] programming languages, has been published in a public GitHub repository under the permissive terms of the MIT License [61]. A conscious effort has been made to keep RKTK easy to build, install, and use, while simultaneously providing maximum performance for long-running computations. The C++ components of RKTK build cleanly on modern GNU and Clang toolchains without external dependencies, while the Julia components provide a natural interface and quick JIT compilation times for interactive use. It is our hope that the free availability of high-performance software will encourage a new generation of mathematicians and computational scientists to explore the landscape of Runge–Kutta methods.

## 2.1 The BFGS Algorithm

Our basic approach to solving the overdetermined system of Runge–Kutta order conditions (2.1) is to recast them as a nonlinear optimization problem.

**Definition 2.1.** Let $p, s \in \mathbb{N}$. The $(p, s)$ ***Runge–Kutta residual function*** and $(p, s)$ ***Runge–Kutta error function*** are the functions $R_{p,s}, E_{p,s} : \mathbb{R}^{s \times s} \times \mathbb{R}^s \to \mathbb{R}$ defined by:

$$R_{p,s}(A, \mathbf{b}) := \sum_{k=1}^{p} \sum_{\tau \in T_k} \left( \mathbf{b} \cdot \mathbf{\Phi}_A(\tau) - \frac{1}{\tau!} \right)^2 \tag{2.2}$$

$$E_{p,s}(A, \mathbf{b}) := \sum_{\tau \in T_{p+1}} \frac{1}{\sigma(\tau)^2} \left( \mathbf{b} \cdot \mathbf{\Phi}_A(\tau) - \frac{1}{\tau!} \right)^2 \tag{2.3}$$

In these terms, an $s$-stage Runge–Kutta method of order $p$ is nothing more than a root of $R_{p,s}$. Since the distinction between $A$ and $\mathbf{b}$ is unimportant from this perspective, we will gather their entries into a single vector $\mathbf{x} \coloneqq (A, \mathbf{b}) \in \mathbb{R}^{s(s+1)/2}$ consisting of the lower-triangular entries of $A$, listed in row-major order, followed by the entries of $\mathbf{b}$. Note that the omission of the factor of $1/\sigma(\tau)^2$ makes $R_{p,s}$ slightly faster to numerically evaluate.

Because $R_{p,s}$ is a non-negative function, all of its roots are local minima which can be searched for using standard nonlinear optimization algorithms. Our method of choice is the BFGS algorithm [9, 28, 29, 57], a well-known quasi-Newton algorithm that constructs a finite-difference approximation of the inverse Hessian matrix $H_k \approx (\nabla^2 R_{p,s}(\mathbf{x}_k))^{-1}$ by measuring differences in its gradient $\nabla R_{p,s}(\mathbf{x}_k)$ at successive points in the search. At each step, we perform a symmetric rank-two update of the form $H_{k+1} \coloneqq H_k + \mathbf{u}\mathbf{u}^T - \mathbf{v}\mathbf{v}^T$ that enforces the secant condition $H_{k+1}(\mathbf{x}_{k+1} - \mathbf{x}_k) = \nabla R_{p,s}(\mathbf{x}_{k+1}) - \nabla R_{p,s}(\mathbf{x}_k)$ while minimizing $\|H_{k+1} - H_k\|$ in a certain weighted Frobenius norm.

Our choice of the BFGS algorithm is informed by several factors:

- Direct application of Newton's method, either to $R_{p,s}$ or to the system of Runge–Kutta order conditions (2.1), would be unsuitable because we expect the minima of $R_{p,s}$ to occur in positive-dimensional loci. Hence, the Hessian matrix $\nabla^2 R_{p,s}$ is expected to be rank-deficient.

- Regularized Newton iteration, as performed in [41], is a possible fix, but the numerical evaluation of $\nabla^2 R_{p,s}$ and solution of the resulting system of linear equations was found to be too costly to be worthwhile.

- We have empirically found that the BFGS update formula often produces reasonable step directions for objective functions with degenerate Hessian matrices, significantly outperforming gradient descent. We typically observe a linear rate of convergence to a local minimum, with quadratic convergence sometimes occurring within a small neighborhood.

- It is straightforward to construct a high-performance implementation of the BFGS algorithm that uses cache-friendly memory access patterns and avoids unnecessary dynamic memory allocation. These optimizations yield a significant speedup on modern CPUs with multi-level cache hierarchies.

RKTK's implementation of the BFGS algorithm uses a quadratic line search procedure that guarantees a reduction in the value of $R_{p,s}$ at each step. It also includes a slight modification of the BFGS strategy. During each step, the quadratic line search procedure is run twice: once using the BFGS search direction $-H_k \nabla R_{p,x}(\mathbf{x}_k)$, and once using the gradient descent search direction $-\nabla R_{p,x}(\mathbf{x}_k)$. If the gradient search direction outperforms the BFGS search direction (i.e., achieves a larger reduction in the value of $R_{p,s}$), then we set $H_{k+1}$ to identity matrix. Otherwise, we update $H_{k+1}$ as usual. The idea is that a long-running BFGS search can occasionally jump between regions of space where the curvature of the objective function changes drastically. When this happens, the approximate Hessian mapped out from the old region can be misleading in the new region, so we throw out the old information and start fresh.

## 2.2 Parallelization and Data Layout

The numerical evaluation of the objective function $R_{p,s}$ is a basic task that will constantly be performed in any nonlinear optimization algorithm. This involves the calculation of all Butcher weights $\mathbf{\Phi}_A(\tau)$ for all rooted trees $\tau$ of order $\leq p$. The recursive structure of Butcher weights allows for significant optimizations to this procedure, providing order-of-magnitude speedups over a naïve implementation.

The first observation is that calculating $\mathbf{\Phi}_A(\tau)$ for all rooted trees $\tau$ involves a great deal of repeated work which can be reused. For example, the matrix-vector product $\mathbf{c} := A\mathbf{1}$ can be used to evaluate the elementwise vector-vector product $\mathbf{c} \odot \mathbf{c}$, which itself can be reused in $A(\mathbf{c} \odot \mathbf{c})$. Indeed, every Butcher weight $\mathbf{\Phi}_A(\tau)$ corresponding to a rooted tree $\tau$ can be written as a matrix-vector product $\mathbf{\Phi}_A(\tau) = A\mathbf{\Phi}_A(\tau)$ if $\tau = [\tau']$ is a one-legged tree, or a vector-vector product $\mathbf{\Phi}_A(\tau) = \mathbf{\Phi}_A([\tau_1']) \odot \mathbf{\Phi}_A([\tau_2', \dots, \tau_n'])$ if $\tau = [\tau_1', \tau_2', \dots, \tau_n']$ is a multi-legged tree. This leads to a natural dependency graph structure on the set of rooted trees of order $\leq p$:

$\cdots \longrightarrow$ 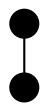 $(\mathbf{c}) \longrightarrow$ 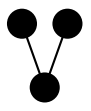 $(\mathbf{c} \odot \mathbf{c}) \longrightarrow$ 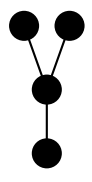 $(A(\mathbf{c} \odot \mathbf{c})) \longrightarrow \cdots$ (2.4)

By performing a breadth-first scan of this dependency graph, we can partition it into a collection of $p$ slices such that:

- The first slice consists solely of the trivial weight $\mathbf{\Phi}_A(\bullet) = \mathbf{1}$.
- For $k > 1$, every Butcher weight in slice $k$ can be evaluated in one step from the Butcher weights in slices $1, 2, \dots, k-1$.

This slicing scheme confers two major benefits. First, it fixes a pattern of data reuse,

ensuring that no matrix-vector or vector-vector product is ever needlessly evaluated more than once. Second, it allows the evaluation of $R_{p,s}$ to be parallelized in each slice, since the breadth-first construction prevents data dependencies between Butcher weights in the same slice. Because each weight has been written as a single matrix-vector or vector-vector product, these tasks can easily be distributed to threads in a load-balanced fashion, providing excellent parallel efficiency with nearly-linear speedups on multi-core systems.

The next observation is that the strictly lower-triangular structure imposed on $A$ forces most Butcher weights to have identically zero entries. For example, the Butcher weight $A^n\mathbf{1}$, which corresponds to a rooted path on $n$ vertices, consists of $n$ zero entries followed by $s - n$ nonzero entries. Because the structures of rooted trees are known ahead of time, we can arrange for the evaluation of $R_{p,s}$ to only store the nonzero entries of Butcher weights explicitly. This provides considerable savings in both memory footprint and computational workload (by omitting needless multiplications by zero).

## 2.3 Dual Number Arithmetic

The use of a quasi-Newton nonlinear optimization algorithm requires us to repeatedly evaluate the gradient $\nabla R_{p,s}$ of the Runge–Kutta residual function. However, the recursive structure of the Butcher weights $\mathbf{\Phi}_A(\tau)$ makes it difficult to construct an explicit formula for their derivatives. Fortunately, the *dual numbers* provide a simple and efficient computational tool for evaluating derivatives without the use of explicit formulae or finite-difference approximations. This construction can be regarded as a simple form of forward-mode automatic differentiation [30, 31, 32, 53].

**Definition 2.2.** The ***dual numbers*** $\mathbb{D} \coloneqq \mathbb{R}[\delta]/\langle \delta^2 \rangle$ are a 2-dimensional commutative algebra over $\mathbb{R}$ obtained by adjoining a formal element[1] $\delta$ having the property that $\delta^2 = 0$. In analogy with the complex numbers $\mathbb{C} \coloneqq \mathbb{R}[i]/\langle i^2 + 1 \rangle$, we call $\delta$ the ***dual unit***, and we say that a dual number $z = a + b\delta$ has ***real part*** $a$ and ***dual part*** $b$.

The utility of the dual numbers is that their arithmetic simulates differentiation for polynomial functions, as codified in the following proposition.

**Proposition 2.3.** Let $n \in \mathbb{N}$ and $p \in \mathbb{R}[x_1, \dots, x_n] \subseteq \mathbb{D}[x_1, \dots, x_n]$. The induced map $p : \mathbb{D}^n \to \mathbb{D}$ satisfies $p(\mathbf{x} + \mathbf{y}\delta) = p(\mathbf{x}) + (\nabla p(\mathbf{x}) \cdot \mathbf{y})\delta$ for all $\mathbf{x}, \mathbf{y} \in \mathbb{R}^n$.

[1] The dual unit $\delta$ is conventionally denoted by $\epsilon$ elsewhere in the literature, but I will use $\epsilon$ to denote the machine epsilon of a finite-precision floating-point arithmetic system in the next section.

*Proof.* We proceed by induction on the number of variables. The base case $n = 0$, where $p$ is a constant function, holds trivially. For $n > 0$, write $p \in \mathbb{R}[x_1, \dots, x_{n-1}][x_n]$ as a polynomial in $x_n$ whose coefficients are polynomials in $x_1, \dots, x_{n-1}$.

$$p(x_1, \dots, x_n) = \sum_{m=0}^{k} p_m(x_1, \dots, x_{n-1}) x_n^m \tag{2.5}$$

Let $\mathbf{x}, \mathbf{y} \in \mathbb{R}^n$ be given, and let $\mathbf{x}', \mathbf{y}' \in \mathbb{R}^{n-1}$ denote $\mathbf{x}$ and $\mathbf{y}$ with their last component dropped. Then, by the inductive hypothesis,

$$\begin{aligned}
p(\mathbf{x} + \mathbf{y}\delta) &= \sum_{m=0}^{k} p_m(\mathbf{x}' + \mathbf{y}'\delta)(x_n + y_n\delta)^m \\
&= \sum_{m=0}^{k} [p_m(\mathbf{x}') + (\nabla p_m(\mathbf{x}') \cdot \mathbf{y}')\delta](x_n + y_n\delta)^m \\
&= \sum_{m=0}^{k} [p_m(\mathbf{x}') + (\nabla p_m(\mathbf{x}') \cdot \mathbf{y}')\delta](x_n^m + m x_n^{m-1} y_n \delta) \\
&= \sum_{m=0}^{k} p_m(\mathbf{x}') x_n^m + \left[ \sum_{m=0}^{k} p_m(\mathbf{x}') m x_n^{m-1} y_n + \sum_{m=0}^{k} (\nabla p_m(\mathbf{x}') \cdot \mathbf{y}') x_n^m \right] \delta \\
&= p(\mathbf{x}) + \left( \frac{\partial p(\mathbf{x})}{\partial x_n} y_n + [\nabla p(\mathbf{x})]' \cdot \mathbf{y}' \right) \delta \\
&= p(\mathbf{x}) + (\nabla p(\mathbf{x}) \cdot \mathbf{y})\delta.
\end{aligned}$$

In the second-to-last line, $[\nabla p(\mathbf{x})]'$ denotes the gradient of $p$ at $\mathbf{x}$ with its last component dropped. □

In a sense, dual number arithmetic can be regarded as a formalization of the non-rigorous manipulations with infinitesimals that are sometimes seen in introductory calculus and physics classes. The multiplication rule $(a + b\delta)(c + d\delta) = ac + (ad + bc)\delta$ is analogous to the product rule $\mathrm{d}(xy) = (x + \mathrm{d}x)(y + \mathrm{d}y) = xy + x\,\mathrm{d}y + y\,\mathrm{d}x$, where the identity $\delta^2 = 0$ codifies the prescription that "terms of order $\mathrm{d}x\,\mathrm{d}y$ are negligible to first order and can be discarded." Alternatively, $\mathbb{D}$ can be regarded as an algebraic structure on the tangent bundle $T\mathbb{R} \cong \mathbb{R} \times \mathbb{R}$, so that dual numbers are thought of as "real numbers with tangents." Arithmetic on $\mathbb{D}$ emerges from taking the differentials of the usual arithmetic operations $+, \times : \mathbb{R} \times \mathbb{R} \to \mathbb{R}$.

These observations provide conceptual motivation for Proposition 2.3, but its true utility lies in its ability to differentiate computer programs. Suppose we are handed a program $P(\mathbf{x})$ that computes a real-valued polynomial function of $\mathbf{x}$. If we can modify $P$ to accept dual numbers as input, then we can construct a program $\nabla P(\mathbf{x})$ to evaluate its gradient by successively evaluating $P(\mathbf{x} + \mathbf{e}_i\delta)$ for $i = 1, 2, \ldots, n$, where $\mathbf{e}_i$ denotes a standard basis vector. This modification is typically easy in programming languages that support function and operator overloading, since it merely involves a substitution of data types (real $\to$ dual) and the addition of an outer loop over $i = 1, 2, \ldots, n$.

Crucially, the only program instructions that need to be modified in this process are arithmetic instructions. All data and control flow constructs remain untouched. Hence, if $P$ has been written in a vectorized/parallelized/cache-optimized fashion, then $\nabla P$ automatically inherits those characteristics. Moreover, $\nabla P$ is typically not as susceptible to catastrophic cancellation as finite-difference approximations would be.

RKTK optimizes the evaluation of gradients via dual arithmetic by separating the computation of the real and dual parts of $R_{p,s}$ and $E_{p,s}$. This allows us to elide needless recomputation of the real part $R_{p,s}(\mathbf{x})$ when repeatedly evaluating $R_{p,s}(\mathbf{x}+\mathbf{e}_i\delta)$ to obtain the components of the gradient vector.

## 2.4 Extended-Precision Machine Arithmetic

The evaluation of the Runge–Kutta residual function $R_{p,s}$ can experience a significant degree of numerical instability because it involves taking the difference $\mathbf{b} \cdot \mathbf{\Phi}_A(\tau) - 1/\tau!$ of two quantities which are expected to be nearly equal. Unfortunately, the definition $R_{p,s}$ does not appear to admit a numerically stable reformulation which avoids this catastrophic cancellation. Thus, extended-precision floating-point arithmetic is necessary to accurately evaluate $R_{p,s}$.

There exist several widely-used libraries for arbitrary-precision floating-point arithmetic, including the GNU GMP/MPFR libraries and Fredrik Johansson's Arb library. These use integer arithmetic to simulate floating-point operations on abstract numeric data types defined by the library. This approach is effective for extremely high-precision workloads (requiring over 1000 digits), but evaluating $R_{p,s}$ accurately enough for the purposes of numerical search only requires slightly higher precision (e.g., 128–256 bits) than the native 64-bit machine arithmetic provided by most computer systems [38]. At this level, the cost of a external function call[2] into a dynamically-loaded library begins to exceed the cost of the actual work that function performs.

[2]On modern computer systems, an external function call itself is not very expensive. The real cost referred to here is that an external function call creates a hard boundary across which a compiler cannot perform optimizations such as inlining and data-flow analysis. Because a dynamically-linked library can be swapped out at run-time, any function provided by such a library must be treated as a black box that can overwrite register values and modify arbitrary regions of memory.

Fortunately, there exists an alternate approach to extended-precision arithmetic that sidesteps these issues while taking full advantage of the arithmetic features of modern CPUs: namely, the “double-double” arithmetic first introduced by Theodorus Dekker in 1971 [21]. The idea of double-double arithmetic is to represent a single 128-bit floating-point number as an ordered pair $(a, b)$ of *non-overlapping* 64-bit[3] floating-point numbers (i.e., the least significant mantissa bit of $a$ is strictly more significant than than the most significant mantissa bit of $b$). Thus, we effectively work with formal sums of the form $a + b\epsilon$, where $\epsilon = 2^{-53}$ denotes the machine epsilon of 64-bit floating-point arithmetic.

Similarly to dual number arithmetic, operations on double-double numbers are performed by discarding terms proportional to $\epsilon^2$. However, care must be taken to ensure that round-off errors are propagated between components. For example, the sum of two double-double numbers is defined by

$$(a + b\epsilon) + (c + d\epsilon) \coloneqq a \oplus c + (b \oplus d \oplus [(a + c) - (a \oplus c)])\epsilon \tag{2.6}$$

where $a \oplus c$ denotes the floating-point sum of $a$ and $c$, and $(a + c) - (a \oplus c)$ denotes the difference between their exact and floating-point sums, i.e., the round-off error incurred in $a \oplus c$. Special algorithms for the exact computation of round-off errors of this type, presented with proofs, can be found in references [21, 50, 51] and [36]. The latter paper [36] also extends this idea to “quad-double” arithmetic, which instead works with formal sums of the form $a + b\epsilon + c\epsilon^2 + d\epsilon^3$ and discards terms of order $\epsilon^4$ and higher. This provides twice the accuracy (256 bits) at the cost of requiring much more complicated procedures for the correct propagation of round-off errors between adjacent components.

The effectiveness of double-double and quad-double arithmetic is significantly boosted by the availability of *vector instructions* on modern CPUs, which are capable of performing arithmetic operations on several pairs of numbers at once. The algorithms for double-double and quad-double arithmetic are particularly suitable candidates for vectorization, since they involve fixed sequences of arithmetic instructions with no branching or looping control flow. Thus, they can be applied uniformly to an arbitrary number of pairs of operands in parallel.

RKTK further extends these ideas, providing algorithms for extended-precision arithmetic using non-overlapping expansions containing anywhere between 2 and 8 terms. This provides a wide range of possible levels of numerical precision (128–512 bits), allowing for fine-grained control of the precision-time trade-off which can even be adjusted dynamically as a computation proceeds.

[3] The 64-bit binary floating-point format defined by the IEEE 754 standard [38] is often called “double precision.” Thus, a number in double-double representation consists of a pair of “doubles.”

## 2.5 Optimization Strategy

The search for a 16-stage Runge–Kutta method of order 10 was performed in several phases at varying levels of numerical precision. First, an initial batch of approximately 800 candidate methods were optimized as far as possible using 64-bit machine arithmetic. Then, the 8 most promising candidate methods (i.e., with the smallest values of $R_{10,16}$) were selected for further optimization using the extended-precision arithmetic techniques described in Section 2.3. After one week of continuous BFGS iteration, one thread reported a nearly-vanishing residual ($R_{10,16} \approx 10^{-200}$) and output the method presented in Appendix A.1.

The principal error coefficients of this method were then minimized via constrained gradient descent applied to $E_{10,16}$ subject to $R_{10,16} = 0$. After three days of continuous optimization, the method presented in Appendix A.2 was obtained. We note that this method is not, in fact, a local minimum of $E_{10,16}$ on the locally 11-dimensional solution variety defined by $R_{10,16} = 0$, so further improvement is possible. Optimization was halted after three days because at this time, each gradient descent step provided only a marginal decrease in $E_{10,16}$ (roughly one part in 10,000,000). We also attempted to use a constrained BFGS-type quasi-Newton algorithm to minimize $E_{10,16}$, but this method performed consistently worse than constrained gradient descent.

All calculations described above were performed on a standard desktop computer equipped with an AMD Ryzen 7 1700 8-core CPU using less than 500MB of RAM.

# Chapter 3

# Results and Conclusions

In this chapter, we discuss the significance of the Runge–Kutta methods discovered using the techniques presented in Chapter 2. We compare these methods to the previously best-known Runge–Kutta methods of order 10 and suggest potentially fruitful directions for future research.

## 3.1 Comparison to Other Runge–Kutta Methods

We compared the performance of ten different Runge–Kutta methods on the following test problem proposed by Erwin Fehlberg [26].

$$y'(t) = -2ty(t)\log z(t) \qquad y(0) = e \tag{3.1}$$

$$z'(t) = 2tz(t)\log y(t) \qquad z(0) = 1 \tag{3.2}$$

This is a system of two non-stiff differential equations that admits a particularly simple analytic solution:

$$y(t) = \exp\big(\cos(t^2)\big) \tag{3.3}$$

$$z(t) = \exp\big(\sin(t^2)\big) \tag{3.4}$$

Each method was run using a fixed step size from $t = 0$ to $t = 5$. A full list of methods tested is provided in the following table. Note that Euler10 and Euler11 are $10^{\text{th}}$-order and $11^{\text{th}}$-order Richardson extrapolations of Euler's method following the construction described in Proposition 1.28.

| Label | Order | Stages | $\sqrt{E_{p,s}}$ | Reference |
|---|---|---|---|---|
| Zhang10 | 10 | 16 | $1.433 \times 10^{-06}$ | Appendix A.2 |
| Hairer10 | 10 | 17 | $5.271 \times 10^{-06}$ | [33] |
| Feagin10 | 10 | 17 | $2.189 \times 10^{-05}$ | [24] |
| Feagin14 | 14 | 35 | $1.060 \times 10^{-05}$ | [25] |
| RK4 | 4 | 4 | $1.450 \times 10^{-02}$ | Equations (1.11)–(1.15) |
| Euler10 | 10 | 46 | $5.753 \times 10^{-08}$ | Proposition 1.28 |
| Euler11 | 11 | 56 | $5.523 \times 10^{-09}$ | Proposition 1.28 |
| RKCK5 | 5 | 6 | $9.483 \times 10^{-04}$ | [17] |
| DOPRI5 | 5 | 7 | $3.991 \times 10^{-04}$ | [23] |
| RKF8 | 8 | 13 | $1.094 \times 10^{-05}$ | [26] |

A log-log plot of the performance of each method as a function of computational workload on Fehlberg's test problem is shown below. Performance was measured as the number of correct digits in the final result, i.e., $-\log_{10}(\|\mathbf{y}_1 - \mathbf{y}(5)\|)$. Computational workload was measured as the number of times the derivative function $\mathbf{f}(t, \mathbf{y}(t))$ was evaluated by each method.

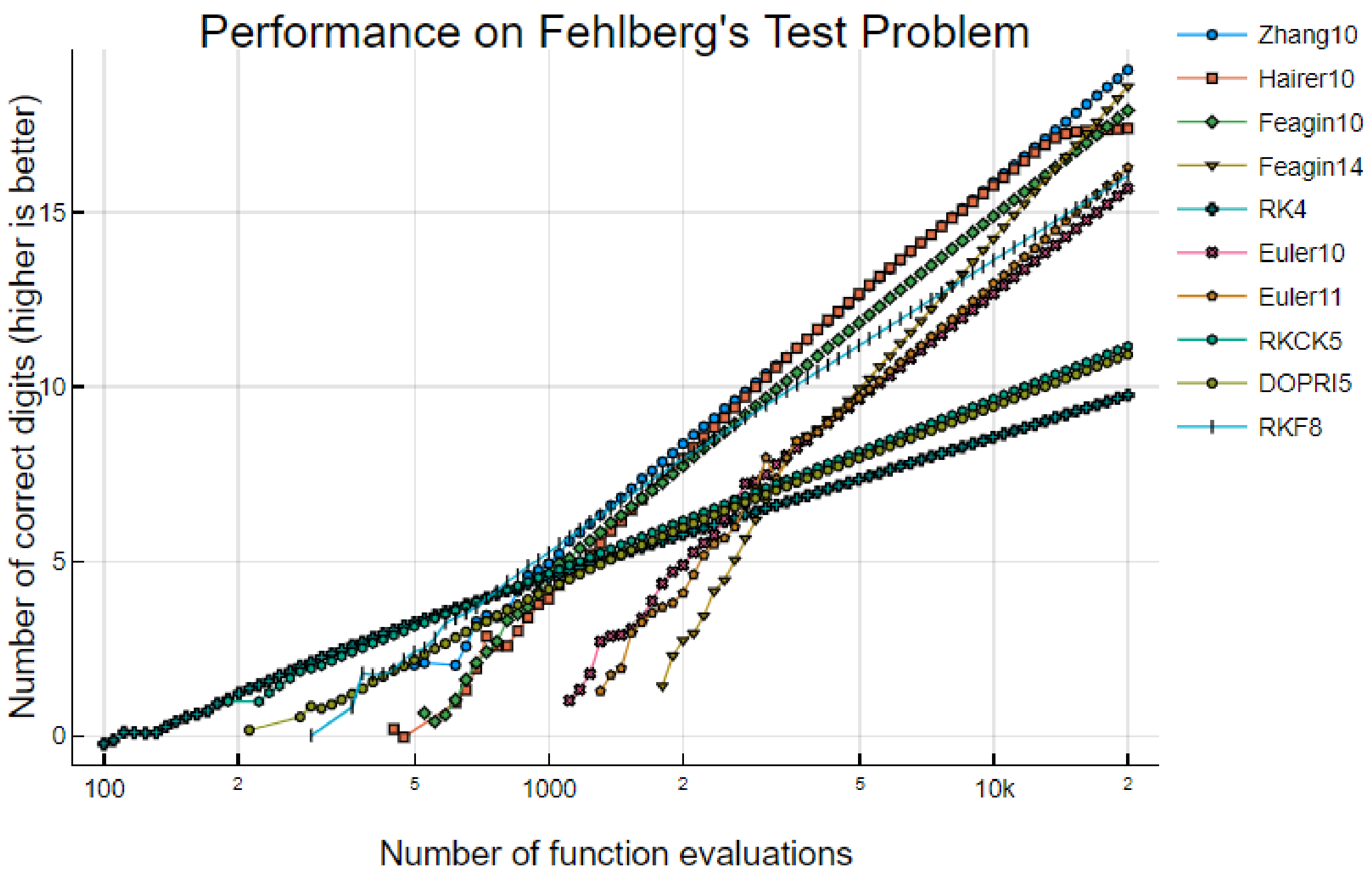

A detailed view of the upper-right region of this plot is provided below.

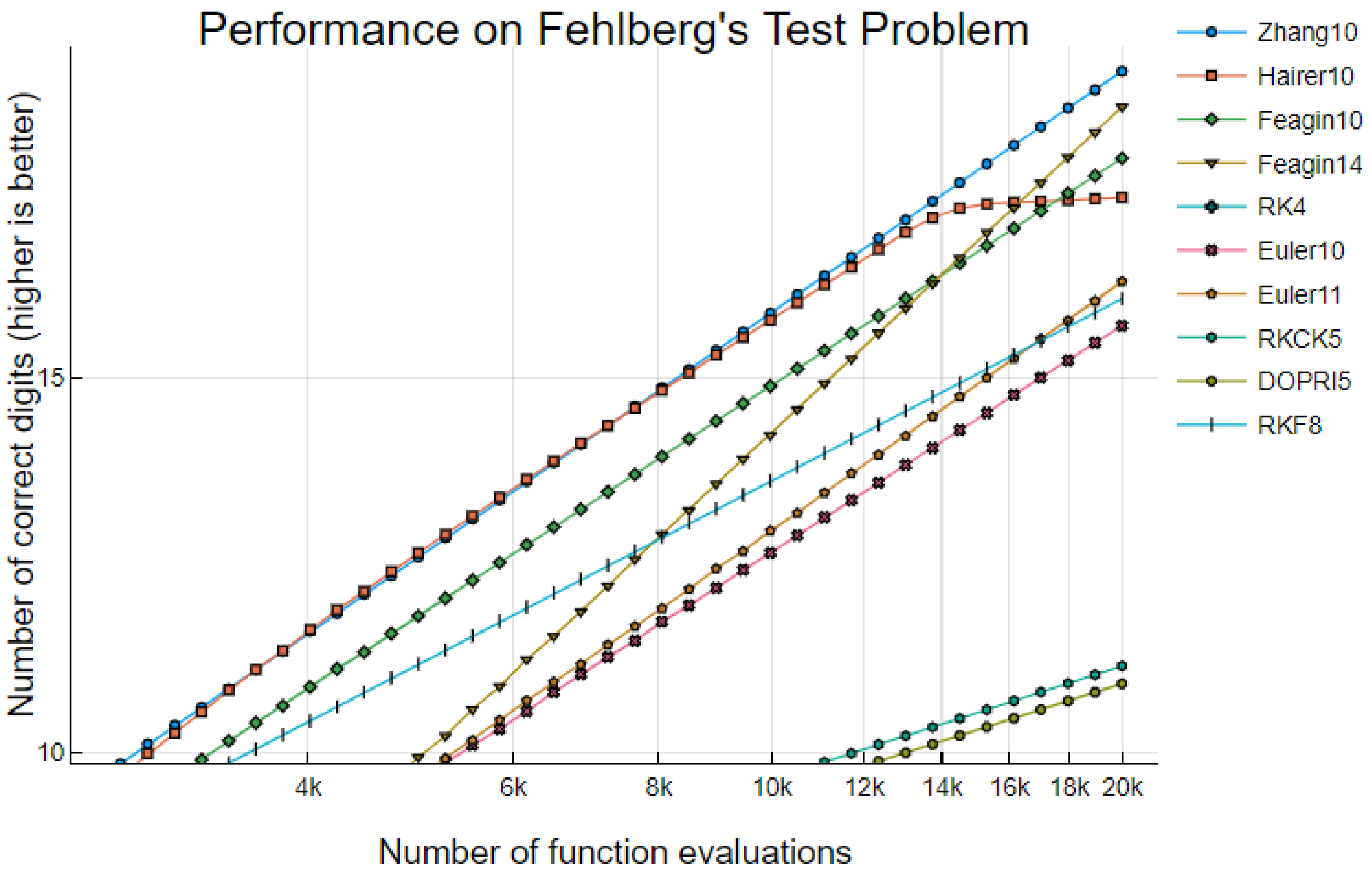

In this regime (≤ 20-digit precision), our optimized method performs comparably to Ernst Hairer's 17-stage method of order 10 and outperforms all others, even Terry Feagin's 35-stage method of order 14 (the highest-order method currently known). These experiments demonstrate that our method is not only a theoretical curiosity, but is suitable for practical computation at 64-bit machine precision.

To conclude this section, we direct attention to a particularly curious feature of the Runge–Kutta methods presented in Appendix A. Although these methods were derived without the assistance of any simplifying assumptions, they nonetheless satisfy Butcher's classical simplifying assumptions $C(2)$ and $D(1)$; namely, $b_2 = 0$ and $c_s = 1$. To date, no Runge–Kutta methods of order $p \geq 6$ are known to violate either of these conditions [41]. It is presently unclear whether these are, in fact, necessary conditions implied by the Runge–Kutta order conditions (1.63). However, we point out that RKTK, having neither of these simplifying assumptions built into its optimization routines, serves as an excellent tool to search for potential counterexamples.

## 3.2 Conclusions and Future Work

In this paper, we presented the numerical discovery of an 11-parameter family of 16-stage Runge–Kutta methods of order 10. We gave an expository account of the theory of Runge–Kutta methods and showed how the structure of this theory can be fruitfully exploited to construct high-performance computational tools for method discovery, analysis, and optimization. The existence of these methods alone is a nontrivial result that constructively improves an upper bound which has not budged in over 40 years [33]. In addition, the discovery of our methods without the help of simplifying assumptions is a novel feature demonstrating that the use of these assumptions can be avoided by making efficient use of modern computer hardware. Finally, we showed that our methods are suitable for practical computation, exhibiting performance that is competitive with previously best-known methods.

We close by stating several natural directions for future research in this area, many of which can benefit from the application of the tools and techniques developed in this paper.

- Continue running RKTK's optimization programs to discover more methods, aiming to further close the existing gaps between upper and lower bounds.
- Use RKTK to optimize the principal error coefficients of existing Runge–Kutta methods and explore their solution manifolds.
- Extend RKTK to derive and optimize methods with other desirable practical properties, including sparsity of the matrix $A$, low-error embedded methods, and the FSAL property.
- Develop tools for the analysis of other families of Runge–Kutta-type methods, including implicit Runge–Kutta methods, Runge–Kutta-Nyström methods, and general linear methods.

Appendix A

# Supplemental Data

## A.1 Coefficients of the first $10^{\text{th}}$-order 16-stage Runge–Kutta method found by RKTK

The following coefficients are also available in plain text format in the RKTK GitHub repository [61].

| | |
|---|---|
| $a_{2,1}$ | +0.0219165081323316356157810083449691695714144018271489194258208394511 |
| $a_{3,1}$ | −3.1742646241245797556990663695838198839944704666223785253697831732740 |
| $a_{3,2}$ | +3.5698359545545936519634451459882514067557621405269776623093968842887 |
| $a_{4,1}$ | −0.0121486945920447046366319234203780234184755658879422852129936963512 |
| $a_{4,2}$ | +0.0353457124363402979980490534278250354091310673029369662254658985799 |
| $a_{4,3}$ | −0.0012805097119639577456361216624778424192410995878457615866513627775 |
| $a_{5,1}$ | +2.5177669651090823650470921704079186390851893568918608559149006983727 |
| $a_{5,2}$ | −1.2469472873778344242622665414211709910996445910577862821396786813361 |
| $a_{5,3}$ | +0.7437014583735782586069289906174905393811638256949786164956076867162 |
| $a_{5,4}$ | −1.3248514658394580730686465923746130892495221906439558596765698273353 |
| $a_{6,1}$ | −1.3970898366224982818785785540247678670536641795425790461393081928548 |
| $a_{6,2}$ | +0.1790892612549676839239583657875869299329769455115140791779835855907 |
| $a_{6,3}$ | +0.1204519669290338645811468189572903513756164407983430228183671149675 |
| $a_{6,4}$ | +1.5385895748357027960127183182921664274348312851017789181232738230645 |
| $a_{6,5}$ | +0.4432689677615069488561490182397866765040326677237410010126698066909 |
| $a_{7,1}$ | −0.11050720319255520608930107538499936804322258913284749994994110280069 |
| $a_{7,2}$ | +0.0220109629323927901678804070829439176979097106256336536474057145352 |
| $a_{7,3}$ | −0.0009868797032498061177712161911126153793978691408133178574298675981 |
| $a_{7,4}$ | +0.1891003283086871655161799706417618033828587962126272130903346128512 |
| $a_{7,5}$ | +0.0028926336333200003964868669155195982701525034214296953610392895360 |
| $a_{7,6}$ | −0.0012490262197379670415555136346612864400551714627422067056880308331 |
| $a_{8,1}$ | +0.7728535804051034554561869781548814589433672657518494111098758828394 |
| $a_{8,2}$ | −0.1091003072167388802201931398434432450891707810070153423324005218399 |
| $a_{8,3}$ | +0.0060429828966800913308068069741553837172780520883279720512774895961 |
| $a_{8,4}$ | −0.9373012883003927999747885222232326172355532324967180199947136534811 |
| $a_{8,5}$ | −0.0105763110905825497882534592343766531807688489485217548490267453366 |
| $a_{8,6}$ | +0.0048733460650133231766196861240605658397231738350894404683511987936 |
| $a_{8,7}$ | +0.5105900100142107500627464720925783230715976644319527725755459006076 |
| $a_{9,1}$ | −0.2166573236838437367028576346134957989229784300131356601226582620726 |
| $a_{9,2}$ | +0.0327877074291082693219833044172952732372818335107349057632275780636 |
| $a_{9,3}$ | +0.0668633811041174873921443220197747470764136318558211968857454833668 |
| $a_{9,4}$ | +0.2816853700756989536400860301736537676690735092070104765996809304924 |
| $a_{9,5}$ | +0.2275458054146548300763203953008141624891606358022088749754055700095 |
| $a_{9,6}$ | −0.0510561826228016381170107965968078350056996933406056062657500397528 |
| $a_{9,7}$ | +0.0216539919920255750555165452317362048273608734521171654023577986193 |
| $a_{9,8}$ | +0.2973205568456127972281357214112773087144383705026457104650026605286 |
| $a_{10,1}$ | $-0.2862290169804059679745396275487759580590817708682509611345472258266$ |
| $a_{10,2}$ | $+0.2922637435448687512458302187759821513933165537934729545118529536639$ |
| $a_{10,3}$ | $-0.0278765549216865687880219359909581877887861448734185906986092653334$ |
| $a_{10,4}$ | $+0.2105422149491281431867787091556708659628718339980781026204138969807$ |
| $a_{10,5}$ | $-0.0842590500040167251101942652654720719130506426744265636892490995298$ |
| $a_{10,6}$ | $+0.0073234684736620187941580672265969854860531811943258301176383157367$ |
| $a_{10,7}$ | $-0.1894463769485046215637979579067624933375201483090174152741998657338$ |
| $a_{10,8}$ | $+0.3863901317935886606528954570382423587345991277645850651912248020540$ |
| $a_{10,9}$ | $+0.0868627705233802058212701109199078722828896838792507152950891990031$ |
| $a_{11,1}$ | $-0.4848260312763151153161850717382625342600713508640649095620254029068$ |
| $a_{11,2}$ | $+0.0039898888119002330673988264460237746471101247247661859036750695341$ |
| $a_{11,3}$ | $+0.0099830137385587367278696358384284842785131676775268531152176540073$ |
| $a_{11,4}$ | $+0.7224924301237864020724462673501658269050504364511758313300935458799$ |
| $a_{11,5}$ | $+0.0638737937556250101220344726040387378955329941630539020482603464971$ |
| $a_{11,6}$ | $-0.0125802839881642401526383102592132575436168374146103807972332217664$ |
| $a_{11,7}$ | $-0.1099693811642481153672225769893826415221813422715396828704681206748$ |
| $a_{11,8}$ | $+0.20970290660290953188923746001096544777714001911191879888059069097646$ |
| $a_{11,9}$ | $-0.0392977633161971682804057462863361012654451323015417035487269067827$ |
| $a_{11,10}$ | $+0.0322027571421586215018438194280037858492595946206450525149138374624$ |
| $a_{12,1}$ | $-0.9994924285514547575349351350543462889531710408276023951905943642470$ |
| $a_{12,2}$ | $+0.1440428799832632827101255886074002919707303181319275019555059160547$ |
| $a_{12,3}$ | $+0.1935193438533814652683334036329461846047140629224828702551340532371$ |
| $a_{12,4}$ | $+1.2374995123579008698887274875847637312912758491430664701470568281725$ |
| $a_{12,5}$ | $+0.5763431827834449981712245805778716143222105264219894582213442485527$ |
| $a_{12,6}$ | $-0.1225525565890229363787934907357897216239902696521444329747603155873$ |
| $a_{12,7}$ | $-0.4250664759461262297145316028877540293729595330514922011521255384584$ |
| $a_{12,8}$ | $+0.7195974108587167382498619080087971584010702483716397435250852332781$ |
| $a_{12,9}$ | $-0.1583051673115796097726940284544055451825367534166606493958281341326$ |
| $a_{12,10}$ | $-0.4164173966122439874074144001830941416031872302898667215930737197866$ |
| $a_{12,11}$ | $-0.0838559203173060157329845301856706108218615934690348735527581381792$ |
| $a_{13,1}$ | $-1.9449356988905405236638623185307672358125122177803270218466754595696$ |
| $a_{13,2}$ | $+0.2435947026713890489560736441180126636320293188353687107753814903908$ |
| $a_{13,3}$ | $+0.1576961206410022374043726231492805496728405935370128656654109360557$ |
| $a_{13,4}$ | $+2.0927679716195473620491453239888007856977480779104024233493936021978$ |
| $a_{13,5}$ | $+0.6524263120891624282993509101970223319110057635768759727583793394383$ |
| $a_{13,6}$ | $-0.0348694789527951985808126301839759088470295708957676191574135740126$ |
| $a_{13,7}$ | $+0.0333334148714097073356428449039300238415868946279672547989483682138$ |
| $a_{13,8}$ | $-0.1259924068416971870226364907352144694143107533385726462158489704019$ |
| $a_{13,9}$ | $-0.4057071515997700290245493649899338019375037512735406294032115147413$ |
| $a_{13,10}$ | $+0.2110769423319483908873873059919781406502945357892260808269618148869$ |
| $a_{13,11}$ | $-0.2433400812037149345796451239633129745675807145281333769021542432344$ |
| $a_{13,12}$ | $+0.2471420033659951073117777392260862452374363622193431640067364495284$ |
| $a_{14,1}$ | $-0.0092046221561684936788476219859530066336542431639053689772545099458$ |
| $a_{14,2}$ | $+0.9540874918506226995502394663812325943078980100545849192468791398517$ |
| $a_{14,3}$ | $-0.0859187767998488930775592358259635288014557091800668090496650894451$ |
| $a_{14,4}$ | $-0.8885800260428398094665145859457603493017051158867444473531393320976$ |
| $a_{14,5}$ | $-0.4003177216167644612082313772091233040824580269086415776998159916640$ |
| $a_{14,6}$ | $-0.0356039647188433844992795088282467815441163055794578187361029755320$ |
| $a_{14,7}$ | $-0.0169764519838735305137822819432966580238350464121719920000915106509$ |
| $a_{14,8}$ | $+0.3296721875987360401643883760322966452765574020667554197802099487896$ |
| $a_{14,9}$ | $+0.4484142499401486136847990584375370600622818901056170280344118613822$ |
| $a_{14,10}$ | −0.2492581183095646788973628371321914102239497286141944650970363635463 |
| $a_{14,11}$ | +0.4066792683268736773674248303828785350775638003321675396262587280600 |
| $a_{14,12}$ | −0.1571885098682438205446900611544347065965905264645032621788374119101 |
| $a_{14,13}$ | +0.0997663242097799373837945551954564332447552735551599713437972177230 |
| $a_{15,1}$ | −1.1686788744160770413014994051629413764118655436582615227803521444666 |
| $a_{15,2}$ | +2.4187132475930298065306983446868228807391598660897483241065173392731 |
| $a_{15,3}$ | −0.0509124971298587256195105493409376812466036119303978484338382988715 |
| $a_{15,4}$ | −1.1308052535593938435879115383728212384083601827662817106309211066978 |
| $a_{15,5}$ | −0.4093638950301158475949722155021051734635416232145904794696589615830 |
| $a_{15,6}$ | −0.2731399180598763256909890348125478905589356625779302260961081509399 |
| $a_{15,7}$ | +0.0528953109949057146489762958043787694814957839796122204254154513127 |
| $a_{15,8}$ | +0.2052253173112451084700580540742786398348143134403631109441438425253 |
| $a_{15,9}$ | +0.7202935710723142977307075852450277828513012242883244779749232513450 |
| $a_{15,10}$ | +0.8461583787899225711020683819939455496595720302218239911408708805123 |
| $a_{15,11}$ | −0.2803297913349801261227309002922209807319120916364468799163503921612 |
| $a_{15,12}$ | −0.4306157124613094578979931262187206184494985516776422324432548634778 |
| $a_{15,13}$ | +0.3449468590342277992521709167271952379770392277260607641423523316338 |
| $a_{15,14}$ | −0.4488154123740200336546940324249223785113735043797828520241254673897 |
| $a_{16,1}$ | +0.8091169256187795848039532797028565479725826694705025407114059045414 |
| $a_{16,2}$ | −0.1178650843090338263897497821322070016838283241427656920055127272375 |
| $a_{16,3}$ | +0.1694292458371054507280277153406831736940917895790760745938873208158 |
| $a_{16,4}$ | −1.0126011391427320291412958156694872994425072413043714126471758127680 |
| $a_{16,5}$ | −0.0679824427636257260042489944186139433981070552665278338474928136620 |
| $a_{16,6}$ | +1.4752905943351533012806579423705970338357995264055551081683369590030 |
| $a_{16,7}$ | +0.6772304270043187398754697224326687350269768505470201002159864075527 |
| $a_{16,8}$ | −0.2628453594795650250645500147137229486273045574953909886402815931832 |
| $a_{16,9}$ | −0.8103278928689712394706589949673207639607356358180986316147484793187 |
| $a_{16,10}$ | −0.3237446998502856420630977939864046803417973043383970960196200765767 |
| $a_{16,11}$ | +0.0254835222696497910245287488575218626612558712242105999605354034688 |
| $a_{16,12}$ | +0.9706715440501786876455755173124670833363451099802348548503803099819 |
| $a_{16,13}$ | −1.2071932695883131176272063064984524728039489984807026159531156467366 |
| $a_{16,14}$ | +1.5252623572889415872307425405434198107718867930864974995371988475269 |
| $a_{16,15}$ | −0.8499247284016005368281477641740051370407094934468425073097840034078 |
| $b_{1}$ | +0.0301200782032484622007933485076449368418567076608727207687269694049 |
| $b_{2}$ | 0.0000000000000000000000000000000000000000000000000000000000000000000 |
| $b_{3}$ | +0.0136040666354818670502332317818773580780095607148577420958730634533 |
| $b_{4}$ | 0.0000000000000000000000000000000000000000000000000000000000000000000 |
| $b_{5}$ | 0.0000000000000000000000000000000000000000000000000000000000000000000 |
| $b_{6}$ | +0.5515323690416867484058433940793350740508945005145689878132567443936 |
| $b_{7}$ | +0.1517807837623530124319039835601067345004412837154457200017321004380 |
| $b_{8}$ | +0.1063656063911084138389660498351490401108933207454179821530385180935 |
| $b_{9}$ | +0.4379736812108413135520924611399533471074194832101827869627419796986 |
| $b_{10}$ | −0.1222844865653052468864691024191391713593537209338937176037364780171 |
| $b_{11}$ | +0.2726693673644476828602050850155321166233647073066483769782057534509 |
| $b_{12}$ | −0.1746665816851846781498336950297707231882036813522129427918082902852 |
| $b_{13}$ | −0.3690902045979758836923665088267880555861225890348542177431259029294 |
| $b_{14}$ | +0.1150989864212221433978098743016724309232597281870315120455365312465 |
| $b_{15}$ | −0.0453658196133234163803081342837171526831358527358479822560229736490 |
| $b_{16}$ | +0.0322621534313995813711300123381440645806765520017830315755819847014 |

## A.2 Coefficients of a $10^{\text{th}}$-order 16-stage Runge–Kutta method optimized by RKTK

The following coefficients are also available in plain text format in the RKTK GitHub repository [61].

| | |
|---|---|
| $a_{2,1}$ | +0.0688809661218865223067709866163293538131515932269828598068246003348 |
| $a_{3,1}$ | −0.8381052035336423753518636620280483869410696889210049308107911986162 |
| $a_{3,2}$ | +1.2536900487169546534492343625929021093721584525964292193606803899247 |
| $a_{4,1}$ | −0.0049094867593268017536963446767904327119018515976089319794041235444 |
| $a_{4,2}$ | +0.0823217303076802157114760304414350071649805075863234687291824351733 |
| $a_{4,3}$ | −0.0085312774264668916510086991483152206399270627617316769429537112940 |
| $a_{5,1}$ | +1.0489319537643595381175146095280468151986564112704572858568112390543 |
| $a_{5,2}$ | −0.7538281773175958124190428032291159017925091699096055719513133160066 |
| $a_{5,3}$ | +0.8052281597406491178047665275688443791431671372879154302150032279268 |
| $a_{5,4}$ | −0.3844632207423856106399215916676820395162389087233623819851545897780 |
| $a_{6,1}$ | −0.2399238343332999527001850107148760502033661458712779930987928130498 |
| $a_{6,2}$ | −0.0636426116392922957110716464377648349095634917957057972402331843864 |
| $a_{6,3}$ | +0.1996754313519789586428397831744866324437219470525002620495210833013 |
| $a_{6,4}$ | +0.6103537950954700215898120246782514690635842960161835994977753010037 |
| $a_{6,5}$ | +0.3785994722499904842640375306432949346214156621503379783849655840043 |
| $a_{7,1}$ | +0.0177883394631617255006766711329999264402040825614508358046806424081 |
| $a_{7,2}$ | −0.0110502190550182547215209343047994002174669982412639837764878832829 |
| $a_{7,3}$ | −0.0043934285505289292963334267463533961891177993647324802301543783691 |
| $a_{7,4}$ | +0.1059752729050901873150565839539349530636382890950153865977110368279 |
| $a_{7,5}$ | +0.0040508906963833073247292990371413035267553284560273690993258200987 |
| $a_{7,6}$ | −0.0016792970913422346102200182541182790808468230583857295791312425004 |
| $a_{8,1}$ | +0.2356604622541853792569269035682497313512758582912074979917610008137 |
| $a_{8,2}$ | +0.0889336268970155965611293086879433570148072406867365744055899691756 |
| $a_{8,3}$ | +0.0413883470987685854516951834146572205141793251842480702603992216410 |
| $a_{8,4}$ | −0.8529030360326306991230961409359935924003175037568143989392658084562 |
| $a_{8,5}$ | −0.0230808754811362556300837986487013527727297763333937307585941452949 |
| $a_{8,6}$ | +0.0096859218859185251459080774595498997953921146500417887845940353067 |
| $a_{8,7}$ | +0.8014497793419689386483379470104767969914199068020516966819646048920 |
| $a_{9,1}$ | +0.0904869376070533068179502765319509084236256152123573721170994960752 |
| $a_{9,2}$ | +0.0320442720247922624237800179739017119198480429527935894722150528562 |
| $a_{9,3}$ | +0.1243376821589105628009505331482701499224744691370015296271736685262 |
| $a_{9,4}$ | −0.3073152175503815810494623201158308497790095186166091648004244461852 |
| $a_{9,5}$ | +0.1644784368116026881016316155162742744858144412077187323701085938216 |
| $a_{9,6}$ | −0.0410068116734476843964655729510615092248565822457165056270171514712 |
| $a_{9,7}$ | +0.4066198993077319845697715354546192410285117624781133929144085449008 |
| $a_{9,8}$ | +0.1866998712463193918001877760787186903409643438289378766317615332030 |
| $a_{10,1}$ | −0.1281584913772105802837912268411472995759545228211031009288904085574 |
| $a_{10,2}$ | −0.0949429224253224554374888160619277034677176082796718500059750768045 |
| $a_{10,3}$ | −0.1445034451182625989491704572990043918630408164231444547699078628206 |
| $a_{10,4}$ | +0.9213494070440691236360615681920579870691526036309646476310606945908 |
| $a_{10,5}$ | −0.1326530105494904632048446798250796817262014576110728469185109217348 |
| $a_{10,6}$ | +0.0166148587263146434831901258689786099156477131870690008787188615196 |
| $a_{10,7}$ | −0.6446177324399182556882653865950215617579255679072138748923829666015 |
| $a_{10,8}$ | +0.47148463683422024098060807066053294761028203774917809473770514989982 |
| $a_{10,9}$ | +0.15101154448891262356107150246546481622684638215041867281807172181887 |
| $a_{11,1}$ | −0.35394262889926481874318272132638227200275604000107564690424492655525 |
| $a_{11,2}$ | −0.12999250061102426001102206164524015480707041252098090235235458563233 |
| $a_{11,3}$ | +0.07296717474130959358123820183761890201581557542115580687780801792533 |
| $a_{11,4}$ | +1.23409288203437673148811178545199717478084359792767549031487434785244 |
| $a_{11,5}$ | +0.24325479314414790350916184694938182359002523242070889454561579643422 |
| $a_{11,6}$ | −0.04886415340347670065787549059652108014832442895629947389530738182699 |
| $a_{11,7}$ | −0.51442371407978914516038594645191171303742863244306956194414424856855 |
| $a_{11,8}$ | +0.07088500449199304304537010242765865457286791828303557337358385517366 |
| $a_{11,9}$ | −0.20555333994396552677698205301698855743095009068049025897065388713688 |
| $a_{11,10}$ | +0.04716132770900545782293703693524094489806604422476436750471220364000 |
| $a_{12,1}$ | −1.28883596401187309154174375754588391722337985378291577953064993833099 |
| $a_{12,2}$ | −0.28020869761589692854933312826845810393637896090137224607343344880155 |
| $a_{12,3}$ | +2.95502531060581466025398140893693575414694538631086569687435523799644 |
| $a_{12,4}$ | +2.68729452804277601048934070442766694554427130448420986662047690713133 |
| $a_{12,5}$ | +4.94062342991935573838476417426809163760371481941041460313006887242755 |
| $a_{12,6}$ | −1.10669748003212502578007105340045335120159611347391447923604138591655 |
| $a_{12,7}$ | −1.01499184284449721513261150169983810954496707583685170059698090784088 |
| $a_{12,8}$ | +2.53750982609373480524394613209930567128221625331336116097749979112888 |
| $a_{12,9}$ | −3.05085435525303312587778339284072228780737127406865428592625671405455 |
| $a_{12,10}$ | −4.48248369012964780236263511347668602256466957962418723733767273268600 |
| $a_{12,11}$ | −1.21257843188510861623323738621796488386816954981559478788233918659322 |
| $a_{13,1}$ | −0.26440344286774078027313975052115341819194404740705505988980388884477 |
| $a_{13,2}$ | −0.06681118811399606825394685412417487195704333031683770523718698572711 |
| $a_{13,3}$ | +0.22983178806827308655694675660046084355033982772309698316158066683966 |
| $a_{13,4}$ | +0.64074149645736181719485655167408177784648088775776468738158085006288 |
| $a_{13,5}$ | +0.42942911081576828545053652299580217770451868606900245607955202491955 |
| $a_{13,6}$ | −0.00967078876662303036856558223558993132212111489612543886857914521711 |
| $a_{13,7}$ | +0.01613666672999681392501812041282804272240296713357286534780333216133 |
| $a_{13,8}$ | −0.03596314318487727365094853269870982256844289210631058204493077377822 |
| $a_{13,9}$ | −0.04902670091896965612083428201101205091198148478269054303875029017011 |
| $a_{13,10}$ | +0.01388537673498896790224810676002758506392838309962714417674884945700 |
| $a_{13,11}$ | −0.02329475109524907760636494276820464664967810607545806961371703044088 |
| $a_{13,12}$ | +0.00411071472206196739547606137318514638902786727009559929394539277822 |
| $a_{14,1}$ | +1.28995927881082034516023610876593719136951515200841564716720963602988 |
| $a_{14,2}$ | +0.13503265202676397192358271319486347590315093450300003507414087911922 |
| $a_{14,3}$ | −1.57533735338626298015011591165511659112640458688783385930436332029911 |
| $a_{14,4}$ | −1.13614144125767380945133294734170682018114470434355261134291967031477 |
| $a_{14,5}$ | −2.60629976327340284860734294234508518773200013764798796731373995419355 |
| $a_{14,6}$ | −0.14225514857221282929776869245957057127882456911938845124276329363255 |
| $a_{14,7}$ | −1.30050279397530175137190231972314623324238865480815606289913301764255 |
| $a_{14,8}$ | +2.70797163114645615543070920761013389431867941711011998785845670786144 |
| $a_{14,9}$ | +2.43554296946750176014394493560037156798756615586273835416220476865955 |
| $a_{14,10}$ | −0.26458182516446999987007272710827096690076304619434648405874330108399 |
| $a_{14,11}$ | +0.28431747275350668428305940521170347967298804085514279367418757904855 |
| $a_{14,12}$ | −0.01567336501605444936187349334932741547111318613954509638698908027766 |
| $a_{14,13}$ | +0.60355253162364202926624736416406789911182794847681800316234125803399 |
| $a_{15,1}$ | −0.82170145239485721248182614307298014363501501247665665859542963907666 |
| $a_{15,2}$ | +1.25591221619821840505806060191670198289896021578755350597122093077688 |
| $a_{15,3}$ | −0.01615284116525821921698589364201132298146948312403205059324494512955 |
| $a_{15,4}$ | −0.02297038987201048533939731478067813642773684981394747854111883904 86 |
| $a_{15,5}$ | −0.0279603437670711148257496960961580980595745773093270356124706066005 |
| $a_{15,6}$ | −0.7509429872708069228856785296294584378472839219042646801373266261127 |
| $a_{15,7}$ | −0.0083102908927800860354453379896957289733706035082138452126825386851 |
| $a_{15,8}$ | +0.0278552386207705781939523176878812789961835461187166573061992447059 |
| $a_{15,9}$ | +0.0402142153933365150373043270398897702435872154842688675351860125371 |
| $a_{15,10}$ | +1.6171364636009542589380411657231441232825488795528975532635259510494 |
| $a_{15,11}$ | −1.3962837774124321431762360366303216121808468002871487754219706220971 |
| $a_{15,12}$ | −0.0142443915906106420830484872694615503917980253911513028166554176376 |
| $a_{15,13}$ | +0.7570905433664326516997037491793619675729164415810244981375129030418 |
| $a_{15,14}$ | −0.2240573576305733047853240218713603700660122610342949667328566164147 |
| $a_{16,1}$ | +0.2772640185318553107109525084420866295476920984386065478227043958611 |
| $a_{16,2}$ | +0.0945381837072822961931214529086851801995457726805579043890542226962 |
| $a_{16,3}$ | +0.8417275429545428654198078897506589717982597224943613268427410367692 |
| $a_{16,4}$ | −0.9066525983284447594329845771395468398064904333093933561404894346537 |
| $a_{16,5}$ | −0.0933485803392326372033922293829223154382675486838667779309117805586 |
| $a_{16,6}$ | +4.0888775891411401096557453879105412055178576584197616506162448374040 |
| $a_{16,7}$ | +0.7953998649942899068707949254598231473949581278257343653554622082283 |
| $a_{16,8}$ | −0.0485917514587563625121870285439552941803211435327462064100785971137 |
| $a_{16,9}$ | +0.1481985145749843088714306457597076940224003500401848660964944810970 |
| $a_{16,10}$ | −1.7702114806218959705675512380830383857698297803456203334302766048702 |
| $a_{16,11}$ | +1.8382108449200817473203066689344636236906057571221332890210541429828 |
| $a_{16,12}$ | +0.0490934446343091672217990928284453150845428120596792346162257342099 |
| $a_{16,13}$ | −3.8037882306700753807124346741068027860479259679550378174513093079786 |
| $a_{16,14}$ | +0.3315723987233985414920503523421767302965629675490758524826797556329 |
| $a_{16,15}$ | −0.8422897607634791433274591770803228763095903928034305458795950897068 |
| $b_1$ | +0.0318192745802340975941941994408892683959596745880488931384186585135 |
| $b_2$ | 0.0000000000000000000000000000000000000000000000000000000000000000000 |
| $b_3$ | +0.0468136928901842195439860702517242419186583685612007227797627359428 |
| $b_4$ | 0.0000000000000000000000000000000000000000000000000000000000000000000 |
| $b_5$ | 0.0000000000000000000000000000000000000000000000000000000000000000000 |
| $b_6$ | +1.3755353617054574901312633974759815109485248982304168872322054848629 |
| $b_7$ | +0.1750656714396424894359074039754808659518684031779558910861533882005 |
| $b_8$ | +0.1492479865300845523448905416854461097586215404544416648521849956519 |
| $b_9$ | +0.2718099231266237214235598816258283028246784494991958450761042129677 |
| $b_{10}$ | −0.1707794051136107538719129444061190270676275720177995974191446193595 |
| $b_{11}$ | +0.3003551976884852181945325569434449914376550238040993982294537643519 |
| $b_{12}$ | −0.0101425440320279863640523131141883093066791048716041001130711499952 |
| $b_{13}$ | −1.1917781578209319095661311581106188144083874482181901757635993666342 |
| $b_{14}$ | +0.0363913935417071365309524609730234759365268240889557980207339099299 |
| $b_{15}$ | −0.0468331608651064555672004144382370964301132585070407799722521943159 |
| $b_{16}$ | +0.0324947663292581801700103176973444800403142012103195528530501798837 |